\documentclass[]{article}
\begin{document}
\newtheorem{proposition}{Proposition}[section]
\newtheorem{definition}{Definition}[section]
\newtheorem{lemma}{Lemma}[section]

\title{\bf Representations$^{6-th}$ of Lie Algebras}
\author{Keqin Liu\\Department of Mathematics\\The University of British Columbia\\Vancouver, BC\\
Canada, V6T 1Z2}
\date{Augest, 2011}
\maketitle

\begin{abstract} We introduce representations$^{6-th}$ of Lie algebras, and study the counterparts of the P-B-W Theorem and the Hopf algebra structure for the enveloping algebras of Lie algebras in the context of representations$^{6-th}$ of Lie algebras.


\end{abstract}

\bigskip
Throughout, an associative algebra always means an associative algebra having an identity, a homomorphism from an associative algebra to an associative algebra always preserves the identity, and all vector spaces are finite-dimensional vector spaces over a field $\mathbf{k}$.

\medskip
Let $V$ be a vector space, and let $End(V)$ be the associative algebra of all linear transformations from $V$ to $V$. It is well-known that $End(V)$ is a Lie algebra with respect to the following bracket $[\, ,\,]$
\begin{equation}\label{eq0}
[f, g]:=fg-gf \quad\mbox{for $f$, $g\in End(V)$}.
\end{equation}
This Lie algebra is denoted by $g\ell (V)$. In this paper, the bracket defined by (\ref{eq0}) will be called the {\bf ordinary bracket}, and a homomorphism from a Lie algebra $\mathcal{L}$ to the Lie algebra  $g\ell (V)$ will be called an {\bf ordinary representation} of the Lie algebra $\mathcal{L}$ on $V$. The ordinary representation theory of reductive Lie algebras has been studied successfully for a long time, but the ordinary representation theory of non-reductive Lie algebras is still beyond reach. Based on the fact that a non-reductive complex Lie algebra does not have a faithful ordinary finite-dimensional irreducible representation, we know that if 
$f: \mathcal{L}\to g\ell (V)$ is a faithful ordinary finite-dimensional representation of a non-reductive complex Lie algebra $\mathcal{L}$, then all linear transformations $f(x)$ with $x\in \mathcal{L}$  have a non-zero proper common invariant subspace $W$ of $V$. In other words, the image of any faithful ordinary finite-dimensional representation of a non-reductive complex Lie algebra is inside the subalgebra 
$End _{_W}(V)$  of $End(V)$, where  $W$ is a  non-zero proper subspace of a vector space $V$ and $End _{_W}(V)$ is the set   of all  linear transformations from $V$ to $V$ which have $W$ as an invariant subspace, i.e.,
$$
End _{_W}(V):=\{\, f \,|\, \mbox{$f\in End(V)$ and $f(W)\subseteq W$}\,\}.
$$
Although every non-reductive complex Lie algebra has a faithful ordinary finite-dimensional representation, it is a very difficult problem to construct these faithful ordinary finite-dimensional representations for non-reductive complex Lie algebras. Since it does not appear that the solution to the difficult problem will soon be found, we search for 
the new way of representing Lie algebras faithfully by using linear transformations. 

\medskip
Except the ordinary bracket, there are other brackets which make $End _{_W}(V)$ into a Lie algebra. In fact, if we fix a scalar $k\in\mathbf{k}$ and a linear transformation $q$ satisfying
$$
q(W)=0\quad\mbox{and}\quad\mbox {$q(v)-v\in W$ for all $v\in V$},
$$
then $End _{_W}(V)$ becomes a Lie algebra with respect to the following square bracket 
$[\, ,\,]_{6, k}$
\begin{equation}\label{eq-1}
[f, g]_{6, k}:=fg-gf-fgq+gfq+ kfqg - kgqf,
\end{equation}
where $f$, $g\in End _{_W}(V)$. This Lie algebra is denoted by 
$g\ell _{q, W}^{[ \,, \, ] _{6,k}}(V)$. Replacing $g\ell (V)$ by 
$g\ell _{q, W}^{[ \,, \, ] _{6,k}}(V)$, we get the notion of representations$^{6-th}$ of Lie algebras, which is the new representations of Lie algebras introduced in this paper. Using 
representations$^{6-th}$ of Lie algebras seems to be a good way of representing non-reductive complex Lie algebras with Abelian radicals faithfully by using linear transformations which have a non-zero proper common invariant subspace.

\medskip
The Poincar$\acute{e}$-Birkhoff-Witt Theorem (P-B-W Theorem) and the Hopf algebra structure for the enveloping algebras of Lie algebras are of great importance to the ordinary representations of Lie algebras. To initiate the study of representations$^{6-th}$ of Lie algebras, it is natural to study the counterparts of the P-B-W Theorem and the Hopf algebra structure for the enveloping algebras of Lie algebras in the context of representations$^{6-th}$ of Lie algebras. The purpose of this paper is to present our results in this study.

\medskip
This paper consists of seven sections.
In section 1 and section 2,  based on  a class of associative algebras which are called invariant algebras by us, we introduce the notion of representations$^{6-th}$ of Lie algebras, and give the connections between representations$^{6-th}$ and the ordinary representations of Lie algebras. From section 3 to section 5, we construct the enveloping$^{6-th}$ algebra of a Lie algebra by using free invariant algebras. In section 6, we prove the extended P-B-W Theorem in the context of representations$^{6-th}$ of Lie algebras. The last section is devoted to the discussion of the Hopf-like algebra structures for the enveloping$^{6-th}$ algebras of Lie algebras.

\medskip
\section{Invariant Algebras}

\medskip
We begin this section with the notion of invariant algebras, which comes from our search for the generalizations of Lie algebras. The most important example of invariant algebras is the subalgebra  $End _{_W}(V)$ of the associative algebra $End(V)$. 

\medskip
Let $A$ be an associative algebra with an idempotent $q$. The set
\begin{equation}\label{eq1.1}
(A, q):=\{\, x, \,|\, \mbox{$x\in A$ and $qxq=qx$}\,\}
\end{equation}
is clearly an associative subalgebra of $A$. The associative algebra $(A, q)$ is called the {\bf (right) invariant algebra} induced by the idempotent $q$. One property of invariant algebras is the following

\begin{proposition}\label{pr01.4} If $q$ is an idempotent of an associative algebras $A$, then $(A,q)$ becomes a Lie algebra under the following square bracket
\begin{equation}\label{eq01.79}
[x, y ]_{6, k} := xy-yx-xyq+yxq+kxqy-kyqx,
\end{equation}
where $x$, $y\in (A,q)$, and $k$ is a fixed scalar in the field $\mathbf{k}$.
\end{proposition}

\medskip
\noindent
{\bf Proof} By (\ref{eq1.1}) and (\ref{eq01.79}), we have
\begin{eqnarray}
[[x,y]_{6, k},z]_{6, k}&=&[[x,y],z]-[[x,y],z]q+\nonumber\\
\label{eq01.92}&& +\underbrace{k^2xqyz}_{1(\ref{eq01.93})}-\underbrace{k^2yqxz}_{2(\ref{eq01.94})}-\underbrace{k^2zqxy}_{3(\ref{eq01.94})}+\underbrace{k^2zqyx}_{4(\ref{eq01.93})},
\end{eqnarray}
where $[\, , \,]$ is the ordinary bracket.

\medskip
It follows from (\ref{eq01.92}) that 
\begin{eqnarray}
[[y,z]_{6, k},x]_{6, k}&=&[[y,z],x]-[[y,z],x]q+\nonumber\\
\label{eq01.93}&& +\underbrace{k^2yqzx}_{5(\ref{eq01.94})}-\underbrace{k^2zqyx}_{4(\ref{eq01.92})}-\underbrace{k^2xqyz}_{1(\ref{eq01.92})}+\underbrace{k^2xqzy}_{6(\ref{eq01.94})}
\end{eqnarray}
and
\begin{eqnarray}
[[z,x]_{6, k},y]_{6, k}&=&[[z,x],y]-[[z,x],y]q+\nonumber\\
\label{eq01.94}&& +\underbrace{k^2zqxy}_{3(\ref{eq01.92})}-\underbrace{k^2xqzy}_{6(\ref{eq01.93})}-\underbrace{k^2yqzx}_{5(\ref{eq01.93})}+\underbrace{k^2yqxz}_{2(\ref{eq01.92})}
\end{eqnarray}

Adding (\ref{eq01.92}), (\ref{eq01.93}) and (\ref{eq01.94}), we get that 
$$
[[x,y]_{6, k},z]_{6, k}+[[y,z]_{6, k},x]_{6, k}+[[z,x]_{6, k},y]_{6, k}=0.
$$
This completes the proof of Proposition~\ref{pr01.4}.

\hfill\raisebox{1mm}{\framebox[2mm]{}}

\medskip
We finish this section with the following

\begin{definition}\label{def01.1} Let $(A,q_{_A})$ and $(B,q_{_B})$ be invariant algebras. A linear map $\phi: (A,q_{_A}) \to (B,q_{_B})$ is called an {\bf invariant homomorphism} if
$$
\phi(xy)=\phi(x)\phi(y)\quad\mbox{for $x$, $y\in (A,q_{_A})$}
$$
and
$$ \phi(1_{_A})=1_{_B}, \qquad \phi(q_{_A})=q_{_B},$$
where $1_{_A}$ and $1_{_B}$ are the identities of $A$ and $B$, respectively. A bijective invariant homomorphism is called an {\bf invariant isomorphism}.
\end{definition}

\medskip
\section{The 6-th General Linear Lie Algebra}

We introduce the notion of  representations$^{6-th}$ of Lie algebras by using invariant algebras. A representation$^{6-th}$ of a Lie algebra $L$ on $V$ is roughly a combination of two ordinary representations of the Lie algebra $L$ on $V$. We finish this section with the connection between  
representations$^{6-th}$ and the ordinary representations of Lie algebras.

\medskip
Let $End(V)$ be the associative algebra of all linear transformations from a vector space $V$ to $V$, and let $I$ (or $I_V$) be the identity linear transformation of $V$. If $W$ is a subspace of $V$, then an linear transformation $q$ satisfying
\begin{equation}\label{eq01.2}
q(W)=0\quad\mbox{and}\mbox (q-I)(V)\subseteq W
\end{equation}
is an idempotent because
$$
q^2(v)-q(v)=q(q-I)(v)\in q(W)=0\quad\mbox{for all $v\in V$.}
$$
An element $q$ of $End(V)$ satisfying (\ref{eq01.2}) is called a {\bf $W$-idempotent}. 

It is clear that the invariant algebra $(End(V), q)$ induced by the $W$-idempotent $q$ is the set 
$End _{_W}(V)$ of all linear transformations on $V$ which have a common invariant subspace $W$, i.e.,
\begin{equation}\label{eq01.3}
(End(V), q)=\{\, f \,|\, \mbox{$f\in End(V)$ and $f(W)\subseteq W$}\,\}.
\end{equation}
The invariant algebra $(End(V), q)$ is called the {\bf linear invariant algebra over $V$} induced by the $W$-idempotent $q$. 

By Proposition~\ref{pr01.4}, $(End(V), q)$ becomes a Lie algebra under the following square bracket:
\begin{equation}\label{eq1.4}
[ x, y ] _{6, k} =xy-yx-xyq+yxq+ kxqy - kyqx
\end{equation}
where $x$, $y\in (End(V), q)$ and $k$ is a fix scalar in the field $\mathbf{k}$. The Lie algebra
is denoted by 
$g\ell _{q, W}^{[ \,, \, ] _{6,k}}(V)$ and is called the {\bf $6$-th general linear Lie algebra} induced by $(q, W)$. Having fixed a basis of $V$:
$$
v_1, \dots , v_n, \quad\underbrace{w_1, \dots , w_m}_{\mbox{a basis of $W$}},
$$
we identify an element of $x\in End_{q, W}(V)$ with a lower triangular block matrix:
$$
x=\left(\begin{array}{c|c}A_{n\times n} & 0\\\hline B_{m\times n} &C_{m\times m}  \end{array}\right),
$$
where $B_{m\times n}\in M_{m,n}(\mathbf{k})$ and $M_{m,n}(\mathbf{k})$ is the ordinary associative matrix algebra of $m\times n$ matrices whose entries are in the field $\mathbf{k}$. In particular, there exists a $T_{m\times n}\in M_{m,n}(\mathbf{k})$ such that
\begin{equation}\label{eq1.5}
q=\left(\begin{array}{c|c}I_{n\times n} & 0\\\hline T_{m\times n} &0 \end{array}\right),
\end{equation}
where $I_{n\times n}$ is the $n\times n$ identity matrix

It follows from (\ref{eq1.4}) that 
\begin{equation}\label{eq1.6}
g\ell _{n,m}^{[ T] _{6,k}}  (\mathbf{k}):=
\left\{\, \left.
\left(\begin{array}{c|c}A & 0\\\hline B &C  \end{array}\right)
\,  \right| \, 
\begin{array}{cc}A\in M_{n, n}(\mathbf{k}),  \\
B\in M_{m, n}(\mathbf{k}),\,  C\in M_{m, m}(\mathbf{k}) \end{array}\,\right\}
\end{equation}
is a Lie algebra under the following square bracket:
\begin{equation}\label{eq1.7}
\Bigg [ \left(\begin{array}{c|c}A_x & 0\\\hline B_x &C_x  \end{array}\right), \, \left(\begin{array}{c|c}A_y & 0\\\hline B_y &C_y  \end{array}\right)\Bigg] _{6,k}=
\end{equation}
\begin{displaymath}
\left(\begin{array}{c|c}k[A_x, A_y] & 0\\\hline 
kB_xA_y-kB_yA_x+kC_xTA_y-kC_yTA_x-C_xC_yT+C_yC_xT &[C_x , C_y] \end{array}\right),
\end{displaymath}
where $[A_x, A_y]:=A_xA_y-A_yA_x$ is the ordinary square bracket. This Lie algebra
$\left(g\ell _{n,m}^{[ T] _{6,k}}  (\mathbf{k}), [ \, , \, ] _{6,k} \right)$ is called the 
{\bf 6-th general matrix Lie algebra}. 

\begin{definition}\label{def1.1} Let $W$ be a subspace of a vector space $V$ over a field $\mathbf{k}$ and let $q$ be a $W$-idempotent. A Lie algebra homomorphism $\varphi$ from a Lie algebra $(\mathcal{L},\,[ \,, \,] )$ to the $6$-th general linear Lie algebra 
$g\ell _{q, W}^{[ \,, \, ] _{6,k}}(V)$ is called a {\bf representation$^{6-th}$} of $\mathcal{L}$ on $V$ induced by $(q, W)$.
\end{definition}

\medskip
A representation$^{6-th}$ of a Lie algebra $\mathcal{L}$ is closely related to two ordinary representations of $\mathcal{L}$. In fact,  if $\varphi : \mathcal{L}\to g\ell _{q, W}^{[ \,, \, ] _{6}}(V)$ is a representation$^{6-th}$ of a Lie algebra $\mathcal{L}$, then both 
$f: \mathcal{L}\to g\ell \left(\frac{V}{W}\right)$ and $g: \mathcal{L}\to g\ell \left(W\right)$
are ordinary representations of the Lie algebra $\mathcal{L}$, where $f$ and $g$ are defined by
$f(x)(v+W):=\varphi (x)(v)+W$ and $g(x)(w):=\varphi (x)(w)$ for $x\in \mathcal{L}$, $v\in V$ and $w\in W$. Also, we have the following

\begin{proposition}\label{pr1.1} Let $(\mathcal{L}, [\, , \,])$ be a Lie algebra over a field $\mathbf{k}$.
If $\varphi : L\to g\ell _{q, W}^{[ \,, \, ] _{6, k}}(V)$ is an representation$^{6-th}$ of $\mathcal{L}$ induced by $(W, q)$, then $f$ and $g$ are two ordinary representations of the Lie algebra $\mathcal{L}$ on $V$, where
\begin{equation}\label{eq1.18} 
f(x)=\varphi (x)+kq\varphi (x) -\varphi (x) q \quad\mbox{and}\quad \quad g(x)=k\varphi (x)q
\end{equation}
for $x\in \mathcal{L}$.
\end{proposition}

\medskip
\noindent
{\bf Proof} Note that $q\varphi (x)q=q\varphi (x)$ for $x\in L$. We now prove that  $f$ and $g$ are two ordinary representations of the Lie algebra $\mathcal{L}$ on the vector space $V$.

\medskip
For $x$, $y\in \mathcal{L}$, we have
\begin{eqnarray*}\label{eq1.24} 
&&[f(x), f(y)]\nonumber\\
&=&\varphi (x)\varphi (y)-\varphi (x)\varphi (y) q
+k^2q\varphi (x)\varphi (y)-\varphi (y)\varphi (x)+\nonumber\\
&&\quad +\varphi (y) \varphi (x)q-k^2q\varphi (y)\varphi (x)=f([x, y]),
\end{eqnarray*}
which proves that $f: L\to g\ell (V)$ is an ordinary representations of $\mathcal{L}$.

\medskip
Similarly, we have
\begin{eqnarray*}
&&g([x, y])=k\varphi ([x, y]) q =k[\varphi(x), \varphi(y)]_{6, k}q\\
&=&k(k\varphi (x)q\varphi (y)-k\varphi (y)q\varphi (x))=
k\varphi (x)q\cdot k\varphi (y)q-k\varphi (y)q\cdot k\varphi (x)q\\
&=&g(x)g(y)-g(y)g(x)=[g(x), g(y)],
\end{eqnarray*}
which proves that $g: \mathcal{L}\to g\ell (V)$ is an ordinary representations of $\mathcal{L}$.

\hfill\raisebox{1mm}{\framebox[2mm]{}}

\medskip
\section{Free Invariant Algebras}

\medskip
Let $(\hat{A}, \hat{q})$ be an invariant algebra, and let $\hat{i}$ be a map from a set $X$ to $(\hat{A}, \hat{q})$. The pair $\Big( (\hat{A}, \hat{q}), \hat{i}\Big)$ is called the {\bf free invariant algebra} generated by the set $X$ if the following universal property holds: given any invariant algebra $(A, q)$ and any map $\theta: X\to (A, q)$ there exists a unique invariant homomorphism $\hat{\theta}: (\hat{A}, \hat{q})\to (A, q)$ such that $\hat{\theta} \hat{i}=\theta$; that is, the following digram is commutative
$$
\begin{array}{ccc}
(\hat{A}, \hat{q})&\stackrel{\hat{\theta}}{\longrightarrow}&(A, q)\\
\shortstack{$\hat{i}$}\Bigg\uparrow&&\Bigg\uparrow\shortstack{$\theta$}\\&&\\
X&\stackrel{}{=}&X
\end{array}
$$

The free invariant algebra generated a set is clearly unique. 

\medskip
We now construct free invariant algebras over a field $\mathbf{k}$. Let 
$X:=\{\, x_j \,|\, j\in J \,\}$ be a set, and let $\tilde{q}$ be a symbol which is not an element of $X$. Let $T(V)$ be the tensor algebra based on a vector space $V$, where 
$V=\displaystyle\bigoplus_{_{j\in J}}\mathbf{k}x_j \bigoplus \mathbf{k}\tilde{q}$ is the vector space over $\mathbf{k}$ with a basis $X\bigcup\{\tilde{q}\}$. Recall that the tensor algebra $T(V)$ has the universal property: given any associative algebra $A$ over $\mathbf{k}$ and any $\mathbf{k}$-linear map $\phi: V\to A$, there exists a unique associative algebra homomorphism 
$\hat{\phi}: T(V)\to A$ such that the following digram is commutative
$$
\begin{array}{ccc}
V&\stackrel{}{\hookrightarrow}&T(V)\\
\shortstack{$\phi$}\Bigg\downarrow&&\Bigg\downarrow\shortstack{$\hat{\phi}$}\\&&\\
A&\stackrel{}{=}&A.
\end{array}
$$
Let $I$ be the ideal of $T(V)$ generated by
\begin{equation}\label{eq2.1}
\{\, \tilde{q}\otimes \tilde{q}-\tilde{q}, \,\, \tilde{q}\otimes a\otimes\tilde{q}-\tilde{q}\otimes a \,|\, a\in T(V) \,\}.
\end{equation}
Let $\hat{A}:=\displaystyle\frac{T(V)}{I}$ and $\hat{q}:=\tilde{q}+I$. Then $(\hat{A}, \hat{q})$ is an invariant algebra.

\begin{proposition}\label{pr2.1} The pair $\Big( (\hat{A}, \hat{q}), \hat{i}\Big)$ is the free invariant algebra generated by the set $X=\{\, x_j \,|\, j\in J \,\}$, where the map
$\hat{i}: X\to \Big( (\hat{A}, \hat{q}), \hat{i}\Big)$ is defined by
\begin{equation}\label{eq2.2}
\hat{i}(x_j):=x_j+I \qquad\mbox{for $j\in J$.}
\end{equation}
\end{proposition}

\medskip
\noindent
{\bf Proof} Let $(A, q)$ be an invariant algebra and let $\theta : X\to (A, q)$ be a map. Then $\theta$ can be extended a $\mathbf{k}$-linear map from $V$ to $(A, q)$ such that
\begin{equation}\label{eq2.3}
\theta(\tilde{q})=q.
\end{equation}
By the universal property of the tensor algebra $T(V)$, the  $\mathbf{k}$-linear map 
$\theta : V \to (A, q)$ can be extended to an associative algebra homomorphism 
$\theta ': T(V)\to (A, q)$ satisfying that $\theta '|I=0$. Hence, $\theta ': T(V)\to (A, q)$ induces a homomorphism $\hat{\theta}: \hat{A}=\displaystyle\frac{T(V)}{I}\to (A, q)$ such that
\begin{equation}\label{eq2.4}
\hat{\theta}(\hat{q})=\theta '(\tilde{q})=\theta (\tilde{q})=q
\end{equation}
and
\begin{equation}\label{eq2.5}
(\hat{\theta}\hat{i})(x_j)=\hat{\theta}(x_j+I)=\theta '(x_j)=\theta (x_j) \quad\mbox{for $j\in J$.}
\end{equation}
By (\ref{eq2.4}) and (\ref{eq2.5}), $\hat{\theta}: (\hat{A}, \hat{q}) \to (A, q)$ is an invariant homomorphism such that $\hat{\theta}\hat{i}=\theta$. Since the associative algebra $\hat{A}$ is generated by the set
$$\{\, \hat{1}:=1+I,\,\, \hat{q},\,\, i(x_j)\, |\, j\in J\,\},$$
the invariant homomorphism $\hat{\theta}$ satisfying (\ref{eq2.4}) and (\ref{eq2.5}) is unique. This proves that the pair $\Big( (\hat{A}, \hat{q}), \hat{i}\Big)$ has the university property of the free invariant algebra generated by the set $X$

\hfill\raisebox{1mm}{\framebox[2mm]{}}

\bigskip
Let $\hat{x}_j:=\hat{i}(x_j)=x_j+I$ for $j\in J$. The product of two elements $a$ and $b$ of 
$(\hat{A}, \hat{q})$ will be denoted by $ab$. The next proposition gives a basic property of the free invariant algebra $\Big( (\hat{A}, \hat{q}), i\Big)$ generated by the set $X$.

\begin{proposition}\label{pr2.2} The following subset of  $\Big( (\hat{A}, \hat{q}), \hat{i}\Big)$ 
$$
\hat{S}:=\left\{\, \hat{1},\,\, \hat{q},\,\, \hat{x}_{j_1}\cdots\hat{x}_{j_m},\,\,
\hat{x}_{j_1}\cdots\hat{x}_{j_t}\hat{q}\hat{x}_{j_{t+1}}\cdots\hat{x}_{j_m} \,\left|\, 
\begin{array}{c}\{\,x_{j_1}, \cdots , x_{j_m}\,\}\subseteq X ,\\ m\in \mathcal{Z}_{\ge 1}, \,\, m\geq t\geq 0\, \end{array}\right.\,\right\}
$$
is a $\mathbf{k}$-basis of the vector space $(\hat{A}, \hat{q})$.
\end{proposition}

\medskip
\noindent
{\bf Proof} First, we construct a vector space $W$ with the following $\mathbf{k}$-basis:
$$
\mathring{S}:=\left\{\, \mathring{1},\,\, \mathring{q},\,\, \mathring{x}_{j_1}\cdots\mathring{x}_{j_m},\,\,
\mathring{x}_{j_1}\cdots\mathring{x}_{j_t}\mathring{q}\mathring{x}_{j_{t+1}}\cdots\mathring{x}_{j_m} \,\left|\, 
\begin{array}{c}\{\,x_{j_1}, \cdots , x_{j_1}\,\}\subseteq X ,\\ m\in \mathcal{Z}_{\ge 1}, \,\, m\geq t\geq 0\, \end{array}\right.\,\right\}
$$

Recall that the set
$$
S:=\{\, 1, \,\, y_1\otimes y_2\otimes\cdots \otimes y_n\,|\, \mbox{$y_i\in X\bigcup \tilde{q}$ and $n\in \mathcal{Z}_{\ge 1}$}\,\}
$$
is a $\mathbf{k}$-basis of the tensor algebra 
$T(V)=T\left(\displaystyle\bigoplus_{_{j\in J}}\mathbf{k}x_j \bigoplus \mathbf{k}\tilde{q}\right)$.

\medskip
We now define two degrees of an element in $S$ as follows:
\begin{eqnarray*}
&&deg_{_X}(1):=0, \qquad deg_{_{\tilde{q}}}(1):=0,\\
&&deg_{_X}(y_1\otimes\cdots \otimes y_n):=
\left|\{\, y_i \,|\, \mbox{$y_i\in X$ for $n\geq i\geq 1$}\,\}\right|,\\
&&deg_{_{\tilde{q}}}(y_1\otimes\cdots \otimes y_n):=
\left|\{\, y_i \,|\, \mbox{$y_i=\tilde{q}$ for $n\geq i\geq 1$}\,\}\right|,
\end{eqnarray*}
where $|A|$ means the cardinality of a set $A$. It is clear that
$$
deg_{_X}(y_1\otimes\cdots \otimes y_n)+deg_{_{\tilde{q}}}(y_1\otimes\cdots \otimes y_n)=n.
$$

Using the $\mathbf{k}$-basis $S$ of $T(V)$, we have
$$
T(V)=\bigoplus_{m\ge 0} T_m\quad\mbox{and}\quad T_m=\bigoplus_{t\ge 0} T_{m,t},
$$
where
$$
T_m:=\bigoplus_{\begin{array}{c}deg_{_{\tilde{q}}}(y_1\otimes\cdots \otimes y_n)=m,\\ y_i\in X\bigcup \tilde{q}\end{array}}\mathbf{k}y_1\otimes\cdots\otimes y_n 
$$
$$
T_{m,t}:=\bigoplus_{\begin{array}{c} deg_{_{\tilde{q}}}(y_1\otimes\cdots \otimes y_n)=m,\\ 
deg_{_X}(y_1\otimes\cdots \otimes y_n)=t,\\ y_i\in X\bigcup \tilde{q}\end{array}}
\mathbf{k}y_1\otimes\cdots\otimes y_n.
$$
Thus, we get
$$
T(V)=T_0\bigoplus T_1\bigoplus\left(\bigoplus_{m\ge2, \, t\ge 0}T_{m,t}\right)\qquad\mbox{(as $\mathbf{k}$-vector spaces)}.
$$

\medskip
Next, we define a linear map $\tau: T(V)\to T(V)$ by
$$
\tau\left|(T_0\bigoplus T_1)\right.:=id
$$
and
$$
\tau\left|T_{m,t}\right.: T_{m,t}\to T_{1,t}\qquad\mbox{for $m\ge 2$ and $t\ge 0$},
$$
where $\tau\left|T_{m,t}\right.$ is defined by the rule: for $y_1\otimes\cdots \otimes y_n\in T_{m,t}$ with $m\ge 2$ and $t\ge 0$, $\tau(y_1\otimes\cdots \otimes y_n)\in T_{1, t}$ is the element obtained from $y_1\otimes\cdots \otimes y_n$ by keeping the first $\tilde{q}$ (the most left $\tilde{q}$) and deleting the other $\tilde{q}$'s. For example, we have
$$
\tau(\underbrace{x_{j_1}\otimes \tilde{q}\otimes x_{j_2}\otimes \tilde{q}\otimes x_{j_3}}_{\mbox{It is in $T_{2,2}$}}):=\underbrace{x_{j_1}\otimes \tilde{q}\otimes x_{j_2}\otimes x_{j_3}}_{\mbox{It is in $T_{1,2}$}}
\quad\mbox{for $x_{j_1}$, $x_{j_2}$, $x_{j_3}\in X$}.
$$

\medskip
Finally, using $\tau$, we define a linear map
$\sigma : T(V)\to W$ by
$$
\sigma(1):=\mathring{1}, \qquad \sigma(\tilde{q}):=\mathring{q},
$$
$$
\sigma(\underbrace{x_{j_1}\otimes x_{j_2}\otimes\cdots\otimes x_{j_m}}_{\mbox{It is in $T_0$ and $m\ge 1$}}):
=\mathring{x}_{j_1}\mathring{x}_{j_2}\cdots \mathring{x}_{j_m},
$$
$$
\sigma(\underbrace{x_{j_1}\otimes\cdots x_{j_t}\otimes\tilde{q}\otimes x_{j_{t+1}}\otimes\cdots\otimes 
x_{j_m}}_{\mbox{It is in $T_1$ and $m\ge 1$}}):
=\mathring{x}_{j_1}\cdots\mathring{x}_{j_t}\mathring{q}\mathring{x}_{j_{t+1}}\cdots\mathring{x}_{j_m}
$$
and
$$
\sigma\left|T_{m,t}\right.:=\left(\sigma\left|T_1\right.\right)
\left(\tau\left|T_{m,t}\right.\right)\qquad\mbox{for $m\ge 2$ and $t\ge 0$.}
$$

It is easy to check that $\sigma (I)=0$. Hence, the linear map 
$\sigma : T(V)\to W$ induces a linear map 
$\bar{\sigma}: (\hat{A}, \hat{q})=\displaystyle\frac{T(V)}{I}\to W$ such that
$$
\bar{\sigma}(\hat{1})=\sigma(1)=\mathring{1}, \qquad \bar{\sigma}(\hat{q})=\sigma(\tilde{q}):=\mathring{q},
$$
$$
\bar{\sigma}(\hat{x}_{j_1}\cdots\hat{x}_{j_m})=
\bar{\sigma}(x_{j_1}\otimes\cdots \otimes x_{j_m}+I)=\sigma(x_{j_1}\otimes\cdots \otimes x_{j_m})
=\mathring{x}_{j_1}\cdots\mathring{x}_{j_m}
$$
and
\begin{eqnarray*}
&&\bar{\sigma}(\hat{x}_{j_1}\cdots\hat{x}_{j_t}\hat{q}\hat{x}_{j_{t+1}}\cdots\hat{x}_{j_m})
=\bar{\sigma}(x_{j_1}\otimes\cdots x_{j_t}\otimes\tilde{q}\otimes x_{j_{t+1}}\otimes\cdots x_{j_m}+I)\\
&=&\sigma(x_{j_1}\otimes\cdots x_{j_t}\otimes\tilde{q}\otimes x_{j_{t+1}}\otimes\cdots x_{j_m})
=\mathring{x}_{j_1}\cdots\mathring{x}_{j_t}\mathring{q}\mathring{x}_{j_{t+1}}\cdots\mathring{x}_{j_m};
\end{eqnarray*}
that is, $\bar{\sigma}(\hat{S})=\mathring{S}$. Since $\mathring{S}$ is a $\mathbf{k}$-basis of the vector space $W$, the set $\hat{S}$ is $\mathbf{k}$-linear independent. Since 
$(\hat{A}, \hat{q})=\displaystyle\frac{T(V)}{I}$ is clearly spanned by $\hat{S}$, Proposition~\ref{pr2.2} is true.

\hfill\raisebox{1mm}{\framebox[2mm]{}}

\section{Enveloping$^{6-th}$ Algebras}

Let $(A, q)$ be an invariant algebra. By Proposition~\ref{pr01.4}, 
$\Big((A, q), [\, , \,]_{6,k}\Big)$ is a Lie algebra, where $k$ is a fixed nonzero scalar in $\mathbf{k}$ and $[\, , \,]_{6,k}$ is the  following square bracket 
\begin{equation}\label{eq2.9}
[x , y]_{6,k}=xy-yx-xyq+yxq+kxqy-kyqx\quad\mbox{for $x$, $y\in (A, q)$}.
\end{equation}
This Lie algebra  is denoted by $Lie(A, q)$ or $Lie\Big((A, q), [\, , \,]_{6,k}\Big)$. 
In the remaining part of this paper, we assume that the scalar $k$ is a fixed non-zero scalar in the field $\mathbf{k}$.

\begin{definition}\label{def2.1} Let $(\mathcal{L}, [\,, \,])$ be a Lie algebra (arbitrary dimensionality and characteristic). By a {\bf enveloping$^{6-th}$ algebra} of 
$(\mathcal{L}, [\,, \,])$ we will understand a pair $\Big((\mathcal{U}, \bar{q}), i\Big)$ composed of an invariant algebra $(\mathcal{U}, \bar{q})$ together with a map 
$i: \mathcal{L}\to (\mathcal{U}, \bar{q})$ satisfying the following two conditions:
\begin{description}
\item[(i)]  the map 
$i: (\mathcal{L}, [\,, \,])\to 
Lie\Big((\mathcal{U}, \bar{q}), [\,, \,]_{6, k}\Big)$ is a Lie algebra homomorphism; that is, $i$ is linear and
\begin{eqnarray*}
&&i([ x, y ])=[ i(x), i(y) ]_{6, k}=i(x)i(y)-i(y)i(x)+\\
&&\quad -i(x)i(y)\bar{q}+i(y)i(x)\bar{q}+ki(x)\bar{q}i(y)-ki(y)\bar{q}i(x)
\quad\mbox{for $x$, $y\in\mathcal{L}$,}
\end{eqnarray*}
\item[(ii)] given any invariant algebra $(A, q)$ and any Lie algebra homomorphism 
$f: (\mathcal{L}, [\,, \,])\to Lie\Big((A, q), [\,, \,]_{6,k}\Big)$ there
exists a unique invariant homomorphism 
$f': (\mathcal{U}, \bar{q})\to (A, q)$ such that $f=f'i$; that is, the following diagram is commutative:
$$
\begin{array}{ccc}
(\mathcal{U}, \bar{q})&\stackrel{f'}{\longrightarrow}&(A, q)\\
\shortstack{i}\Bigg\uparrow&&\Bigg\uparrow\shortstack{f}\\&&\\
\mathcal{L}&\stackrel{}{=}&\mathcal{L}.
\end{array}
$$
\end{description}
\end{definition}

Clearly, the enveloping$^{6-th}$ algebra of a Lie algebra $\mathcal{L}$, which is also denoted by 
$(\mathcal{U}^{6-th}(\mathcal{L}), \bar{q})$, is unique up to an invariant isomorphism.

\bigskip
We now construct the enveloping$^{6-th}$ algebra of a Lie algebra 
$(\mathcal{L}, [ \, ,\, ])$ over a field $\mathbf{k}$. 
Let $X=\{\, x_j \,|\, j\in J \,\}$ be a basis of 
$(\mathcal{L}, [\,, \,])$. Let $\Big( (\hat{A}, \hat{q}), \hat{i}\Big)$ be the free invariant algebra generated by the set $X$. By Proposition~\ref{pr2.2},
$\hat{X}:=\{\, \hat{x_j}:=\hat{i}(x_j) \,|\, j\in J \,\}$ is a linearly independent subset of $(\hat{A}, \hat{q})$. Hence, $\hat{i}$ can be extended to an injective linear map 
$\hat{i}: \mathcal{L}\to (\hat{A}, \hat{q})$. Let $\hat{x}:=\hat{i}(x)$ for all 
$x\in \mathcal{L}$. Let $R$ be the ideal of $(\hat{A}, \hat{q})$ which is generated by all the elements of the form
$$
\widehat{[ x, y ]}-\hat{x}\hat{y}+\hat{y}\hat{x}+\hat{x}\hat{y}\hat{q}
-\hat{y}\hat{x}\hat{q}-k\hat{x}\hat{q}\hat{y}+k\hat{y}\hat{q}\hat{x},
$$
where $\hat{x}$, $\hat{y}\in\mathcal{L}:=\hat{i}(\mathcal{L})$. Let 
\begin{equation}\label{eq2.10}
\mathcal{U}:=\displaystyle\frac{\hat{A}}{R}, \quad \bar{q}:=\hat{q}+R.
\end{equation}
Define a map $i: \mathcal{L}\to (\mathcal{U}, \bar{q})$ by
\begin{equation}\label{eq2.11}
i(x): =\hat{x}+R \qquad\mbox{for $x\in \mathcal{L}$}.
\end{equation}

We shall now prove

\begin{proposition}\label{pr2.3} The pair $\Big((\mathcal{U}, \bar{q}), i\Big)$ defined by (\ref{eq2.10}) and (\ref{eq2.11}) is an 
\linebreak enveloping$^{6-th}$ algebra for the Lie algebra 
$(\mathcal{L}, [\,, \,])$
\end{proposition}

\medskip
\noindent
{\bf Proof} Clearly, $\Big((\mathcal{U}, \bar{q})\Big)$
is an invariant algebra induced by the idempotent $\bar{q}$, and $i: (\mathcal{L}, [\, , \,])\to 
Lie\Big((\mathcal{U}, \bar{q}), [\, , \,]_{6, k}\Big)$ is a Lie algebra homomorphism.

\medskip
Given any invariant algebra $(A, q)$ and any Lie algebra homomorphism 
$f: (\mathcal{L}, [\, , \,])\to 
Lie\Big((A, q), [\, , \,]_{6, k}\Big)$. We are going to construct an invariant homomorphism $f': (\mathcal{U}, \bar{q})\to (A, q)$ such that $f=f'i$. By 
Proposition~\ref{pr2.2}, the set
\begin{equation}\label{eq2.12}
\hat{S}:=\left\{\,
\hat{x}_{j_1}\cdots\hat{x}_{j_t}\hat{q}^u\hat{x}_{j_{t+1}}\cdots\hat{x}_{j_m}
 \,\left|\, 
\begin{array}{c}\{\,x_{j_1}, \cdots , x_{j_m}\,\}\subseteq X ,\\ m\in \mathcal{Z}_{\ge 0}, \,\, m\geq t\geq 0, \\ u=0,\, 1\,\end{array}\right.\,\right\}
\end{equation}
is a basis of $(\hat{A}, \hat{q})$, where $\hat{q}^0:=\hat{1}$ is the identity of 
$(\hat{A}, \hat{q})$, and
$$
\hat{x}_{j_1}\cdots\hat{x}_{j_t}\hat{q}^u\hat{x}_{j_{t+1}}\cdots\hat{x}_{j_m}:=\left\{
\begin{array}{ll}\hat{q}^u&\quad\mbox{if $m=0$},\\
\hat{q}^u\hat{x}_{j_{1}}\cdots\hat{x}_{j_m}&\quad\mbox{if $t=0$},\\
\hat{x}_{j_1}\cdots\hat{x}_{j_m}\hat{q}^u&\quad\mbox{if $t=m$}.\end{array}\right.
$$
Hence, we can define a linear map 
$\hat{f}: (\hat{A}, \hat{q})\to (A, q)$ by
\begin{equation}\label{eq2.13}
\hat{f}(\hat{x}_{j_1}\cdots\hat{x}_{j_t}\hat{q}^u\hat{x}_{j_{t+1}}\cdots
\hat{x}_{j_m}):=f(x_{j_1})\cdots f(x_{j_t})q^uf(x_{j_{t+1}})\cdots f(x_{j_m}),
\end{equation}
where $u=0$, $1$ and $q^0:=1$ is the identity of $(A, q)$.

\medskip
It is clear that
$\hat{f}: (\hat{A}, \hat{q})\to (A, q)$ is an invariant homomorphism
and $\hat{f}(R)=0$. Thus, 
$\hat{f}: (\hat{A}, \hat{q})\to (A, q)$ induces an invariant homomorphism 
$f':  \left(\mathcal{U}=\displaystyle\frac{\hat{A}}{R}, \bar{q}:=\hat{q}+R\right)\to (A, q)$
such that
$$
f'(\hat{a}+R)=\hat{f}(\hat{a})\qquad\mbox{for $\hat{a}\in (\hat{A}, \hat{q})$.}
$$
Using the equation above, we have
$$
f'i(x)=f'(\hat{x}+R)=\hat{f}(\hat{x})=f(x)\qquad\mbox{for $x\in\mathcal{L}$,} 
$$
which proves that $f'i=f$.
Since the invariant homomorphism $f'$ satisfying $f'i=f$ is unique, the pair $\Big((\mathcal{U}, \bar{q}), i\Big)$ satisfies the two conditions in 
Definition~\ref{pr2.1}. Hence, Proposition~\ref{pr2.3} holds.

\hfill\raisebox{1mm}{\framebox[2mm]{}}

\medskip
\section{Model Monomials}

Let $X=\{\, x_j \,|\, j\in J \,\}$ be a basis of a Lie algebra 
$(\mathcal{L}, [\, , \, ])$, and let $\Big( (\hat{A}, \hat{q}), \hat{i}\Big)$ be the free invariant algebra generated by the set $X$. By Proposition~\ref{pr2.2}, the set $\hat{S}$ given by (\ref{eq2.12}) is a $\mathbf{k}$-basis of $(\hat{A}, \hat{q})$. An element of $\hat{S}$ is called a {\bf monomial}. The {\bf $X$-degree} and the {\bf $\hat{q}$-degree} of a monomial 
$\hat{x}_{j_1}\cdots\hat{x}_{j_t}\hat{q}^u\hat{x}_{j_{t+1}}\cdots
\hat{x}_{j_m}$ is defined by
$$
deg_{_X}(\hat{x}_{j_1}\cdots\hat{x}_{j_t}\hat{q}^u
\hat{x}_{j_{t+1}}\cdots\hat{x}_{j_m}):=m
$$
and
$$
deg_{\hat{q}}(\hat{x}_{j_1}\cdots\hat{x}_{j_t}\hat{q}^u
\hat{x}_{j_{t+1}}\cdots\hat{x}_{j_m}):=\left\{
\begin{array}{ll}(t, m-t)&\quad\mbox{if $u=1$ and}\\
&\quad\mbox{$m\geq t\geq 0$,}\\
(m+1, -1)&\quad\mbox{if $u=0$}.\end{array}\right.
$$

For $m\in\mathcal{Z}_{\ge 0}$, the set $\{\, (t, m-t)\,|\, m+1\geq t\geq 0 \,\}$ is ordered as follows:
\begin{equation}\label{eq2.21}
(0, m)<(1, m-1)<\cdots <(m, 0)<(m+1, -1).
\end{equation}

We suppose now that the set $J$ of indices is ordered.
The {\bf right index} $ind_r(\hat{x}_{j_1}\cdots\hat{x}_{j_t}\hat{q}^u\hat{x}_{j_{t+1}}\cdots
\hat{x}_{j_m})$ and
the {\bf left index} $ind_{\ell}(\hat{x}_{j_1}\cdots\hat{x}_{j_t}\hat{q}
\hat{x}_{j_{t+1}}\cdots
\hat{x}_{j_m})$ of a monomial $\hat{x}_{j_1}\cdots\hat{x}_{j_t}\hat{q}
\hat{x}_{j_{t+1}}\cdots\hat{x}_{j_m}$ are defined as follows:
$$
ind_{r}(\hat{x}_{j_1}\cdots\hat{x}_{j_t}\hat{q}^u\hat{x}_{j_{t+1}}\cdots
\hat{x}_{j_m}):=\left\{
\begin{array}{ll}\displaystyle\sum_{m\geq k>i\geq 1}\eta_{ik}&\quad\mbox{if $u=0$}\\&\\
\displaystyle\sum_{m\geq k>i\geq t+1}\eta_{ik}&\quad\mbox{if $u=1$},\end{array}\right.
$$
and
$$
ind_{\ell}(\hat{x}_{j_1}\cdots\hat{x}_{j_t}\hat{q}
\hat{x}_{j_{t+1}}\cdots
\hat{x}_{j_m}):=\left\{
\begin{array}{ll}0,&\quad\mbox{if $2\geq t\geq 0$}\\&\\
\displaystyle\sum_{t-1\geq k>i\geq 1}\eta_{ik}&\quad\mbox{if $t\ge 3$},\end{array}\right.
$$
where
$$
\eta_{ik}:=\left\{
\begin{array}{ll}0&\quad\mbox{if $i<k$ and $j_k\geq j_i$,}\\
1&\quad\mbox{if $i<k$ and $j_i> j_k$}\end{array}\right.
$$
and
$$
\displaystyle\sum_{p_{_2}\ge k>i\ge p_{_1}}\eta_{ik}:=0 \quad\mbox{if $1\ge |p_{_1}-p_{_2}|$.}
$$

Let
$$
\tilde{T}:=\left\{\,\begin{array}{c}
\hat{q}\hat{x}_{j_1}\cdots\hat{x}_{j_m},\\
\hat{x}_{i_1}\cdots \hat{x}_{i_t}\hat{x}_{j_0}\hat{q}\hat{x}_{j_1}\cdots\hat{x}_{j_m}, \\
\hat{x}_{j_1}\cdots\hat{x}_{j_m} \\\end{array}
 \,\left|\, 
\begin{array}{c}\{\,x_{i_1}, \cdots , x_{i_t}, x_{j_0}, x_{j_1}, \cdots , x_{j_m}\,\}\subseteq X ,\\
\mbox{$i_t\geq\cdots \geq i_1$,\quad 
$j_m\geq \cdots \geq j_1\geq j_0$},\\
\mbox{$t$, $m\in\mathcal{Z}_{\ge 0}$}
\end{array}\right.\,\right\}.
$$
An element of $\tilde{T}$ is called a {\bf model monomial}.

\begin{proposition}\label{pr2.4} Every element of $(\hat{A}, \hat{q})$ is congruent mod $R$ to a $\mathbf{k}$-linear combination of model monomials.
\end{proposition}

\medskip
\noindent
{\bf Proof} It suffices to prove that
\begin{eqnarray}
&&\mbox{\it every monomial whose $X$-degree equals $m$ is congruent $mod \, R$ to}\nonumber\\
\label{eq2.22}&&\mbox{\it a $\mathbf{k}$-linear combination of model monomials for any $m\in\mathcal{Z}_{\ge 0}$.}
\end{eqnarray}

We order the monomials in $\hat{S}$ by $X$-degrees, for a given $X$-degree by $\hat{q}$-degrees, and for a given $X$-degree and a given $\hat{q}$-degree by left or right indices. Note that
$$
(\hat{A}, \hat{q})=\bigoplus_{m\in\mathcal{Z}_{\ge 0}}\hat{A}_m\quad\mbox{and}\quad
\hat{A}_m=\bigoplus_{t=0}^{m+1}\hat{A}_m^{(t, m-t)},
$$
where $\hat{A}_m$ is the subspace spanned by the monomials whose $X$-degree is $m$, and $\hat{A}_m^{(t, m-t)}$ is the subspace spanned by the monomials whose $X$-degree is $m$ and whose 
$\hat{q}$-degree is $(t, m-t)$. For convenience, if $a\in (\hat{A}, \hat{q})$ and $a$ is a sum of some monomials whose $X$-degrees are $m$, we also say that the $X$-degree of $a$ is $m$.

\medskip
We now begin to prove (\ref{eq2.22}) by induction on $m$. It is clear that (\ref{eq2.22}) holds for $m=0$ or $m=1$ because  a monomial whose $X$-degree is $0$ or $1$ is a model monomial. Assume that
\begin{eqnarray}
&&\mbox{\it (\ref{eq2.22}) holds for any monomial whose $X$-degree}\nonumber\\
\label{eq2.23}&&\mbox{\it is less than $m$ with $m\in\mathcal{Z}_{\ge 2}$.}
\end{eqnarray}
We are going to prove that (\ref{eq2.22}) also holds for any monomial whose $X$-degree is $m$.
First, we prove that
\begin{eqnarray}
&&\mbox{\it (\ref{eq2.22}) holds for any monomial whose $X$-degree is $m$ and}\nonumber\\
\label{eq2.24}&&\mbox{\it whose $\hat{q}$-degree is $(0, m)$.}
\end{eqnarray}

A monomial whose $X$-degree is $m$ and whose $\hat{q}$-degree is $(0, m)$  has the form $\hat{q}\hat{x}_{j_1}\cdots\hat{x}_{j_m}$. We use induction on $n:=ind_r(\hat{q}\hat{x}_{j_1}\cdots\hat{x}_{j_m})$ to prove 
(\ref{eq2.22}). Clearly, (\ref{eq2.24}) holds for $n=0$. Assume that
\begin{eqnarray}
&&\mbox{\it (\ref{eq2.24}) holds for any monomial 
$\hat{q}\hat{x}_{h_1}\cdots\hat{x}_{h_m}$ with}  \nonumber\\
\label{eq2.25}&&\mbox{\it $ind_r(\hat{q}\hat{x}_{h_1}\cdots\hat{x}_{h_m})<n$ and $n\in\mathcal{Z}_{\ge 1}$.}
\end{eqnarray}
Consider any monomial $\hat{q}\hat{x}_{j_1}\cdots\hat{x}_{j_m}$ with
$ind_r(\hat{q}\hat{x}_{j_1}\cdots\hat{x}_{j_m})=n$. Since $n\ge 1$, there exists 
$m\ge s+1>s\ge 1 $ such that $j_s>j_{s+1}$. For $x$, 
$y\in (\mathcal{L}, [\, , \,])$, we have
$$
R\ni \hat{q}\left(\widehat{[ x, y]}-[\hat{x}, \hat{y}]_{6, k}\right)
=\hat{q}\widehat{[ x, y]}-k\hat{q}\hat{x}\hat{y}+k\hat{q}\hat{y}\hat{x}
$$
or
\begin{equation}\label{eq2.26}
\hat{q}\hat{x}\hat{y}\equiv \frac{1}{k}\hat{q}[ x, y]\hat{}\,+
\hat{q}\hat{y}\hat{x}\, (mod \, R),
\end{equation}
where $[ x, y]\hat{}:=\widehat{[ x, y]}$. It follows from (\ref{eq2.26}) that
\begin{eqnarray}
&&\hat{q}\hat{x}_{j_1}\cdots\hat{x}_{j_m}
=\hat{q}\hat{x}_{j_1}\cdots\hat{x}_{j_{s-1}}
\Big(\hat{q}\hat{x}_{j_s}\hat{x}_{j_{s+1}}\Big)\hat{x}_{j_{s+2}}\cdots
\hat{x}_{j_m}\nonumber\\
&\equiv&\frac{1}{k}\underbrace{\hat{q}\hat{x}_{j_1}\cdots\hat{x}_{j_{s-1}}\cdot
[ x_{j_s}, x_{j_{s+1}}]\hat{}\cdot \hat{x}_{j_{s+2}}
\cdots\hat{x}_{j_m}}_{(\ref{eq2.27})_1}+\nonumber\\
\label{eq2.27}&& +\underbrace{\hat{q}\hat{x}_{j_1}\cdots\hat{x}_{j_{s-1}}
\Big(\hat{x}_{j_{s+1}}\hat{x}_{j_s}\Big) \hat{x}_{j_{s+2}}
\cdots\hat{x}_{j_m}}_{(\ref{eq2.27})_2}\,\, (mod \, R).
\end{eqnarray}
Since the $X$-degree of the term $(\ref{eq2.27})_1$ is $m-1< m$, the term $(\ref{eq2.27})_1$ is congruent $mod \, R$ to a $\mathbf{k}$-linear combination of model monomials by $(\ref{eq2.23})$. Since the right index of the term $(\ref{eq2.27})_2$ is $n-1< n$, the term $(\ref{eq2.27})_2$ is congruent $mod \, R$ to a $\mathbf{k}$-linear combination of model monomials by $(\ref{eq2.25})$.
Hence, (\ref{eq2.27}) implies that any monomial 
$\hat{q}\hat{x}_{j_1}\cdots\hat{x}_{j_m}$ whose right index equals $n$ is congruent 
$mod \, R$ to a $\mathbf{k}$-linear combination of model monomials. This completes the poof of (\ref{eq2.24}) by induction on $n$.

By the proof of (\ref{eq2.24}), we have
\begin{equation}\label{eq2.28}
\hat{q}\hat{x}_{i_1}\cdots\hat{x}_{i_w}
\equiv \hat{q}\hat{x}_{j_1}\cdots\hat{x}_{j_w}
\, \left(mod \, \bigoplus_{i=0}^{w-1}\hat{A}_i\right),
\end{equation}
where $m\ge w\ge 1$, $\{\, i_1, \dots , i_w\,\}=\{\, j_1, \dots , j_w\,\}$ and 
$j_w\ge \dots\ge j_1$.

\medskip
Similarly, we have
\begin{eqnarray}
&&\mbox{\it (\ref{eq2.22}) holds for any monomial whose $X$-degree is $m$ and}\nonumber\\
\label{eq2.29}&&\mbox{\it whose $\hat{q}$-degree is either $(1, m-1)$ or 
$(2, m-2)$.}
\end{eqnarray}

\medskip
We now prove that
\begin{eqnarray}
&&\mbox{\it (\ref{eq2.22}) holds for any monomial whose $X$-degree is $m$ and }\nonumber\\
\label{eq2.37}&&\mbox{\it whose $\hat{q}$-degree is $(t, m-t)$ with 
$m\geq t\geq 0$. }
\end{eqnarray}

By (\ref{eq2.24}) and (\ref{eq2.29}), (\ref{eq2.37}) holds for $2\ge t\ge 0$. Assume that
\begin{eqnarray}
&&\mbox{\it (\ref{eq2.22}) holds for any monomial whose $X$-degree is $m$ and whose }\nonumber\\
\label{eq2.38}&&\mbox{\it $\hat{q}$-degree is less than $(t+1, m-t-1)$ with 
$m\geq t+1\geq 3$. }
\end{eqnarray}
We are going to prove that
\begin{eqnarray}
&&\mbox{\it (\ref{eq2.22}) also holds for any monomial whose $X$-degree is $m$ and}\nonumber\\
\label{eq2.39}&&\mbox{\it whose $\hat{q}$-degree is  $(t+1, m-t-1)$ with 
$m\geq t+1\geq 3$. }
\end{eqnarray}

Any monomial whose $X$-degree is $m$ and whose $\hat{q}$-degree is  
$(t+1, m-t-1)$ with $m\geq t+1\geq 3$ has he form
$\hat{x}_{j_1}\cdots\hat{x}_{j_t}
\hat{x}_{j_{t+1}}\hat{q}\hat{x}_{j_{t+2}}\cdots\hat{x}_{j_m}$ with 
$m\ge t+1\ge 3$. We shall prove (\ref{eq2.39}) buy induction on the left index 
$n_{\ell}:=ind_{\ell}(\hat{x}_{j_1}\cdots\hat{x}_{j_t}
\hat{x}_{j_{t+1}}\hat{q}\hat{x}_{j_{t+2}}\cdots\hat{x}_{j_m})$, where 
$m\geq t+1\geq 3$. It is easy to check that (\ref{eq2.39}) holds for any monomial whose $X$-degree is $m$, whose $\hat{q}$-degree is  $(t+1, m-t-1)$ with $m\geq t+1\geq 3$ and 
whose left index is $0$. Assume that
\begin{eqnarray}
&&\mbox{\it (\ref{eq2.39}) holds for any monomial whose $X$-degree is $m$, whose }\nonumber\\
\label{eq2.42}&&\mbox{\it $\hat{q}$-degree is $(t+1, m-t-1)$ with 
$m\geq t+1\geq 3$ and whose}\nonumber\\
&&\mbox{\it left index is less that $n_{\ell}$, where $n_{\ell}\in\mathcal{Z}_{\ge 1}$. }
\end{eqnarray}
If $n_{\ell}:=ind_{\ell}(\hat{x}_{j_1}\cdots\hat{x}_{j_t}
\hat{x}_{j_{t+1}}\hat{q}\hat{x}_{j_{t+2}}\cdots\hat{x}_{j_m})\ge 1$, where $m\geq t+1\geq 3$, then there exists $s$ such that 
$t\geq s+1>s\geq 1$ and $j_s>j_{s+1}$. Using
\begin{equation}\label{eq2.43}
\mbox{$\hat{x}\hat{y}\equiv \hat{y}\hat{x}+[x, y]\,\hat{}+\hat{x}\hat{y}\hat{q}-
\hat{y}\hat{x}\hat{q}-k\hat{x}\hat{q}\hat{y}+k\hat{y}\hat{q}\hat{x}
\,(mod \, R)$ for $x$, $y\in \mathcal{L}$,}
\end{equation}
we have
\begin{eqnarray}\label{eq2.44}
&&\hat{x}_{j_1}\cdots\hat{x}_{j_{s-1}}(\hat{x}_{j_s}\hat{x}_{j_{s+1}})
\hat{x}_{j_{s+2}}\cdots\hat{x}_{j_t}
\hat{x}_{j_{t+1}}\hat{q}\hat{x}_{j_{t+2}}\cdots\hat{x}_{j_m}\nonumber\\
&\equiv&\underbrace{\hat{x}_{j_1}\cdots\hat{x}_{j_{s-1}}(\hat{x}_{j_{s+1}}\hat{x}_{j_s})
\hat{x}_{j_{s+2}}\cdots\hat{x}_{j_t}
\hat{x}_{j_{t+1}}\hat{q}\hat{x}_{j_{t+2}}
\cdots\hat{x}_{j_m}}_{(\ref{eq2.44})_1}+\nonumber\\
&& +\underbrace{\hat{x}_{j_1}\cdots\hat{x}_{j_{s-1}}\cdot 
[\hat{x}_{j_s}, \hat{x}_{j_{s+1}}]\,\hat{}\cdot
\hat{x}_{j_{s+2}}\cdots\hat{x}_{j_t}
\hat{x}_{j_{t+1}}\hat{q}\hat{x}_{j_{t+2}}
\cdots\hat{x}_{j_m}}_{(\ref{eq2.44})_2}+\nonumber\\
&& +\underbrace{\hat{x}_{j_1}\cdots\hat{x}_{j_{s-1}}
(\hat{x}_{j_s}\hat{x}_{j_{s+1}}\hat{q})
\hat{x}_{j_{s+2}}\cdots\hat{x}_{j_t}\hat{x}_{j_{t+1}}
\hat{x}_{j_{t+2}}\cdots\hat{x}_{j_m}}_{(\ref{eq2.44})_3}+\nonumber\\
&& -\underbrace{\hat{x}_{j_1}\cdots\hat{x}_{j_{s-1}}
(\hat{x}_{j_{s+1}}\hat{x}_{j_s}\hat{q})
\hat{x}_{j_{s+2}}\cdots\hat{x}_{j_t}\hat{x}_{j_{t+1}}
\hat{x}_{j_{t+2}}\cdots\hat{x}_{j_m}}_{(\ref{eq2.44})_4}+\nonumber\\
 &&-k\underbrace{\hat{x}_{j_1}\cdots\hat{x}_{j_{s-1}}
(\hat{x}_{j_s}\hat{q}\hat{x}_{j_{s+1}})
\hat{x}_{j_{s+2}}\cdots\hat{x}_{j_t}\hat{x}_{j_{t+1}}
\hat{x}_{j_{t+2}}\cdots\hat{x}_{j_m}}_{(\ref{eq2.44})_5}+\nonumber\\
&&+k\underbrace{\hat{x}_{j_1}\cdots\hat{x}_{j_{s-1}}
(\hat{x}_{j_{s+1}}\hat{q}\hat{x}_{j_s})
\hat{x}_{j_{s+2}}\cdots\hat{x}_{j_t}\hat{x}_{j_{t+1}}
\hat{x}_{j_{t+2}}\cdots\hat{x}_{j_m}}_{(\ref{eq2.44})_6}\nonumber\\
&&\, (mod \, R).
\end{eqnarray}
Since 
$$ ind_{\ell}\left((\ref{eq2.44})_1\right)=n_{\ell}-1<n_{\ell},\quad
deg_X\left((\ref{eq2.44})_2\right)=m-1<m$$
and
$$(t+1, m-t-1)> (s+1, m-s-1)
\geq deg_{\hat{q}}\left((\ref{eq2.44})_i\right)
\quad\mbox{for $6\geq i\geq 3$,}$$
each term on the right hand side of (\ref{eq2.44}) is congruent $mod\, R$ to a $\mathbf{k}$-linear combination of model monomials by (\ref{eq2.42}), (\ref{eq2.23}) and (\ref{eq2.38}). This fact and (\ref{eq2.44}) imply that
\begin{eqnarray}
&&\mbox{\it (\ref{eq2.39}) also holds for any monomial 
whose $X$-degree is $m$, }\nonumber\\
\label{eq2.45}&&\mbox{\it whose $\hat{q}$-degree is $(t+1, m-t-1)$ with 
$m\geq t+1\geq 3$ and}\nonumber\\
&&\mbox{\it whose left index is $n_{\ell}$ 
with $n_{\ell}\in\mathcal{Z}_{\ge 1}$. }
\end{eqnarray}

This completes the proof of (\ref{eq2.39}), which implies that  (\ref{eq2.37}) holds.

\bigskip
Similarly, using induction on the right index, we have
\begin{eqnarray}
&&\mbox{\it (\ref{eq2.22}) holds for any monomial $F$ whose $X$-degree is $m$}\nonumber\\
\label{eq2.46}&&\mbox{\it and whose $\hat{q}$-degree is $(m+1, -1)$.}
\end{eqnarray}

By (\ref{eq2.37}) and (\ref{eq2.46}), we get that if (\ref{eq2.23}) holds, then  (\ref{eq2.22}) holds for monomial whose $X$-degree is $m$. This completes the proof of (\ref{eq2.22}) by induction on $m$.

\hfill\raisebox{1mm}{\framebox[2mm]{}}

\section{The Extended$^{6-th}$ P-B-W Theorem}

\medskip
Throughout this section, $(\mathcal{L}, [\, , \,])$ is a Lie algebra over a field $\mathbf{k}$, $X=\{\, x_j \,|\, j\in J \,\}$ is a basis of 
$(\mathcal{L}, [\, , \,])$, $(\hat{A}, \hat{q})$ is the free invariant algebra generated by the set $X$, $\hat{i}$ is the injective linear map 
$\hat{i}: \mathcal{L}\to (\hat{A}, \hat{q})$ which extends the injective map 
$X\to (\hat{A}, \hat{q})$, and $\Big((\mathcal{U}, \bar{q}), i\Big)$ is the enveloping$^{6-th}$ algebra for the Lie algebra 
$(\mathcal{L}, [\, , \,])$, where
$$\mathcal{U}:=\displaystyle\frac{\hat{A}}{R}, \quad \bar{q}:=\hat{q}+R,$$
the map $i: \mathcal{L}\to (\mathcal{U}, \bar{q})$ is given by
$$i(x): =\hat{x}+R \qquad\mbox{for $\hat{x}:=\hat{i}(x)$ and $x\in \mathcal{L}$},$$
and $R$ is the ideal of $(\hat{A}, \hat{q})$ generated by all the elements of the form
\begin{eqnarray}
&&[ x, y]\,\hat{}\, -\hat{x}\hat{y}+\hat{y}\hat{x}+\hat{x}\hat{y}\hat{q}
-\hat{y}\hat{x}\hat{q}-k\hat{x}\hat{q}\hat{y}+k\hat{y}\hat{q}\hat{x}\nonumber\\
\label{eq3.1}&=&[ x, y] \,\hat{}\, - [ \hat{x}, \hat{y}]_{6,k}\qquad
\mbox{with $x$, $y\in\mathcal{L}$.}
\end{eqnarray}
After the set $J$ of indices is ordered, the $\mathbf{k}$-vector space $\mathcal{U}:=\displaystyle\frac{\hat{A}}{R}$ is spanned by the following set of cosets of model monomials
$$
\hat{T}:=\left\{\begin{array}{c}
\hat{q}\hat{x}_{j_1}\cdots\hat{x}_{j_m}+R,\\
\hat{x}_{i_1}\cdots \hat{x}_{i_t}\hat{x}_{j_0}\hat{q}\hat{x}_{j_1}\cdots\hat{x}_{j_m}+R, \\
\hat{x}_{j_1}\cdots\hat{x}_{j_m}+R \\\end{array}
 \,\left|\, 
\begin{array}{c}x_{i_1}, \cdots , x_{i_t}, x_{j_0}, x_{j_1}, \cdots , x_{j_m}\in X,\\
i_t\geq\cdots \geq i_1,\\ 
j_m\geq \cdots \geq j_1\geq j_0,\\
t, m\in\mathcal{Z}_{\ge 0}
\end{array}\right.\right\}.
$$
by Proposition~\ref{pr2.4}. In this chapter we aim to prove that the spanning set $\hat{T}$ is a $\mathbf{k}$-basis of the enveloping$^{6-th}$ algebra $\Big((\mathcal{U}, \bar{q}), i\Big)$ for the Lie algebra $(\mathcal{L}, [\, , \,])$, which is the extended$^{6-th}$ 
P-B-W Theorem.  We will use four subsections to complete the proof of the extended$^{6-th}$ P-B-W Theorem.

\medskip
\subsection{$(u, v)^{[6-th]}$-generators}

Recall that the ideal $R$ of $(\hat{A}, \hat{q})$ is generated by all the elements of the form 
$[x, y] \,\hat{}\, - [\hat{x}, \hat{y}]_{6,k}$ with $x$, $y\in X$. Every generator $[x, y] \,\hat{}\, - [\hat{x}, \hat{y}]_{6,k}$ of the ideal $R$ is a $\mathbf{k}$-linear combination of the elements of the form
\begin{eqnarray*}
&&[x_{j_s}, x_{j_{s+1}}]\,\hat{}\, 
- [\hat{x}_{j_s}, \hat{x}_{j_{s+1}}]_{6, k}
=[x_{j_s}, x_{j_{s+1}}]\,\hat{}\, -\hat{x}_{j_s} \hat{x}_{j_{s+1}}+\\
&&+\hat{x}_{j_{s+1}}\hat{x}_{j_s}
+\hat{x}_{j_s} \hat{x}_{j_{s+1}}\hat{q}-\hat{x}_{j_{s+1}}\hat{x}_{j_s}\hat{q}
-k\hat{x}_{j_s}\hat{q} \hat{x}_{j_{s+1}}+
k\hat{x}_{j_{s+1}}\hat{q}\hat{x}_{j_s},
\end{eqnarray*}
where $x_{j_s}$, $x_{j_{s+1}}\in X$. This property and Proposition~\ref{pr2.2} imply the following

\begin{proposition}\label{pr3.1} Every element of the ideal $R$ is a $\mathbf{k}$-linear combination of the elements of the form
\begin{eqnarray*}
H_{u,v}:&=&
(\hat{x}_{j_1}\cdots\hat{x}_{j_t}\hat{q}^u\hat{x}_{j_{t+1}}\cdots
\hat{x}_{j_{s-1}})\left([ x_{j_s}, x_{j_{s+1}}] \,\hat{}\, 
- [ \hat{x}_{j_s}, \hat{x}_{j_{s+1}}]_{6,k}\right)\cdot\\
&& \cdot (\hat{x}_{j_{s+2}}\cdots\hat{x}_{j_r}\hat{q}^v\hat{x}_{j_{r+1}}\cdots
\hat{x}_{j_m}),
\end{eqnarray*}
where $\{\, x_{j_1}, \cdots , x_{j_m} \,\}\subseteq X$, $m\ge 2$ and $u$, $v=0$, $1$.
\end{proposition}

\medskip
\noindent
{\bf Proof} Since  $\hat{x}_{j_1}\cdots\hat{x}_{j_t}\hat{q}^u\hat{x}_{j_{t+1}}\cdots
\hat{x}_{j_{s-1}}$ or $\hat{x}_{j_{s+2}}\cdots\hat{x}_{j_r}\hat{q}^v\hat{x}_{j_{r+1}}\cdots
\hat{x}_{j_m}$ can represent any element of the $\mathbf{k}$-basis $\hat{S}$ of the vector space $(\hat{A}, \hat{q})$, this proposition is clear.

\hfill\raisebox{1mm}{\framebox[2mm]{}}

\bigskip
The element $H_{u,v}$ in Proposition~\ref{pr3.1} is called a 
{\bf $(u, v)^{[6-th]}$-generator}. $\hat{A}_0\oplus \hat{A}_1$ does not contain any $(u, v)^{[6-th]}$-generator because the $X$-degrees of some terms in a 
$(u, v)^{[6-th]}$-generator $H_{u,v}$ are $2$. Using the property of the idempotent $\hat{q}$, we can rewrite $H_{u,v}$ as follows:
\begin{eqnarray}
H_{00}&=& \underbrace{\hat{x}_{j_1}\cdots\hat{x}_{j_{s-1}}
\left([ x_{j_s}, x_{j_{s+1}}] \,\hat{}\,\right)
\hat{x}_{j_{s+2}}\cdots\hat{x}_{j_m}}_{H_{00}^{(1)}}+\nonumber\\
&&-\underbrace{\hat{x}_{j_1}\cdots\hat{x}_{j_{s-1}}
\left(\hat{x}_{j_s} \hat{x}_{j_{s+1}}\right)
\hat{x}_{j_{s+2}}\cdots\hat{x}_{j_m}}_{H_{00}^{(2)}}+\nonumber\\
&&+\underbrace{\hat{x}_{j_1}\cdots\hat{x}_{j_{s-1}}
\left(\hat{x}_{j_{s+1}}\hat{x}_{j_s}\right)
\hat{x}_{j_{s+2}}\cdots\hat{x}_{j_m}}_{H_{00}^{(3)}}+\nonumber\\
&&+\underbrace{\hat{x}_{j_1}\cdots\hat{x}_{j_{s-1}}
\left(\hat{x}_{j_s}\hat{x}_{j_{s+1}}\hat{q}\right)
\hat{x}_{j_{s+2}}\cdots\hat{x}_{j_m}}_{H_{00}^{(4)}}+\nonumber\\
&&-\underbrace{\hat{x}_{j_1}\cdots\hat{x}_{j_{s-1}}
\left(\hat{x}_{j_{s+1}}\hat{x}_{j_s}\hat{q}\right)
\hat{x}_{j_{s+2}}\cdots\hat{x}_{j_m}}_{H_{00}^{(5)}}+\nonumber\\
&&-k\underbrace{\hat{x}_{j_1}\cdots\hat{x}_{j_{s-1}}
\left(\hat{x}_{j_s}\hat{q}\hat{x}_{j_{s+1}}\right)
\hat{x}_{j_{s+2}}\cdots\hat{x}_{j_m}}_{H_{00}^{(6)}}+\nonumber\\
\label{eq3.2}&&+k\underbrace{\hat{x}_{j_1}\cdots\hat{x}_{j_{s-1}}
\left(\hat{x}_{j_{s+1}}\hat{q}\hat{x}_{j_s}\right)
\hat{x}_{j_{s+2}}\cdots\hat{x}_{j_m}}_{H_{00}^{(7)}},
\end{eqnarray}
\begin{eqnarray}
H_{10}&=&H_{11}= \underbrace{\hat{x}_{j_1}\cdots\hat{x}_{j_t}\hat{q}
\hat{x}_{j_{t+1}}\cdots\hat{x}_{j_{s-1}}
\left([ x_{j_s}, x_{j_{s+1}}] \,\hat{}\,\right)
\hat{x}_{j_{s+2}}\cdots\hat{x}_{j_m}}_{H_{10}^{(1)}}+\nonumber\\
&&-k\underbrace{\hat{x}_{j_1}\cdots\hat{x}_{j_t}\hat{q}
\hat{x}_{j_{t+1}}\cdots\hat{x}_{j_{s-1}}
\left(\hat{x}_{j_s} \hat{x}_{j_{s+1}}\right)
\hat{x}_{j_{s+2}}\cdots\hat{x}_{j_m}}_{H_{10}^{(2)}}+\nonumber\\
\label{eq3.3}&&+k\underbrace{\hat{x}_{j_1}\cdots\hat{x}_{j_t}\hat{q}
\hat{x}_{j_{t+1}}\cdots\hat{x}_{j_{s-1}}
\left( \hat{x}_{j_{s+1}}\hat{x}_{j_s}\right)
\hat{x}_{j_{s+2}}\cdots\hat{x}_{j_m}}_{H_{10}^{(3)}}
\end{eqnarray}
and
\begin{eqnarray}
H_{01}&=& \underbrace{\hat{x}_{j_1}\cdots\hat{x}_{j_{s-1}}
\left([ x_{j_s}, x_{j_{s+1}}] \,\hat{}\,\right)
\hat{x}_{j_{s+2}}\cdots\hat{x}_{j_r}\hat{q}
\hat{x}_{j_{r+1}}\cdots\hat{x}_{j_m}}_{H_{01}^{(1)}}+\nonumber\\
&& -\underbrace{\hat{x}_{j_1}\cdots\hat{x}_{j_{s-1}}
\left(\hat{x}_{j_s} \hat{x}_{j_{s+1}}\right)
\hat{x}_{j_{s+2}}\cdots\hat{x}_{j_r}\hat{q}
\hat{x}_{j_{r+1}}\cdots\hat{x}_{j_m}}_{H_{01}^{(2)}}+\nonumber\\
&& +\underbrace{\hat{x}_{j_1}\cdots\hat{x}_{j_{s-1}}
\left( \hat{x}_{j_{s+1}}\hat{x}_{j_s}\right)
\hat{x}_{j_{s+2}}\cdots\hat{x}_{j_r}\hat{q}
\hat{x}_{j_{r+1}}\cdots\hat{x}_{j_m}}_{H_{01}^{(3)}}+\nonumber\\
&& +\underbrace{\hat{x}_{j_1}\cdots\hat{x}_{j_{s-1}}
\left( \hat{x}_{j_s}\hat{x}_{j_{s+1}}\hat{q}\right)
\hat{x}_{j_{s+2}}\cdots\hat{x}_{j_r}
\hat{x}_{j_{r+1}}\cdots\hat{x}_{j_m}}_{H_{01}^{(4)}}+\nonumber\\
&&-\underbrace{\hat{x}_{j_1}\cdots\hat{x}_{j_{s-1}}
\left(\hat{x}_{j_{s+1}}\hat{x}_{j_s}\hat{q}\right)
\hat{x}_{j_{s+2}}\cdots\hat{x}_{j_r}
\hat{x}_{j_{r+1}}\cdots\hat{x}_{j_m}}_{H_{01}^{(5)}}+\nonumber\\
&&-k\underbrace{\hat{x}_{j_1}\cdots\hat{x}_{j_{s-1}}
\left(\hat{x}_{j_s}\hat{q} \hat{x}_{j_{s+1}}\right)
\hat{x}_{j_{s+2}}\cdots\hat{x}_{j_r}
\hat{x}_{j_{r+1}}\cdots\hat{x}_{j_m}}_{H_{01}^{(6)}}+\nonumber\\
\label{eq3.4}&&+k\underbrace{\hat{x}_{j_1}\cdots\hat{x}_{j_{s-1}}
\left( \hat{x}_{j_{s+1}}\hat{q}\hat{x}_{j_s}\right)
\hat{x}_{j_{s+2}}\cdots\hat{x}_{j_r}
\hat{x}_{j_{r+1}}\cdots\hat{x}_{j_m}}_{H_{01}^{(7)}}.
\end{eqnarray}

The extended$^{6-th}$ P-B-W theorem is a corollary of the following

\begin{proposition}\label{pr3.2} Let $W$ be a $\mathbf{k}$-vector space with a basis
$$
\mathring{T}:=\left\{\,\begin{array}{c}
\mathring{q}\mathring{x}_{j_1}\cdots\mathring{x}_{j_m},\\
\mathring{x}_{i_1}\cdots \mathring{x}_{i_t}\mathring{x}_{j_0}\mathring{q}\mathring{x}_{j_1}\cdots\mathring{x}_{j_m}, \\
\mathring{x}_{j_1}\cdots\mathring{x}_{j_m} \\\end{array}
 \,\left|\, 
\begin{array}{c}\{\,x_{i_1}, \cdots , x_{i_t}, x_{j_0}, x_{j_1}, \cdots , x_{j_m}\,\}\subseteq X ,\\
\mbox{$i_t\geq\cdots \geq i_1$,\quad 
$j_m\geq \cdots\geq j_1\ge j_0$},\\
\mbox{$t$, $m\in\mathcal{Z}_{\ge 0}$}
\end{array}\right.\,\right\},
$$
where $\mathring{x}_{j_1}\cdots\mathring{x}_{j_m}:=\mathring{1}$ for $m=0$, 
$\mathring{x}_{i_1}\cdots \mathring{x}_{i_t}\mathring{x}_{j_0}:=\mathring{x}_{j_0}$ for $t=0$ and $\mathring{q}\mathring{1}:=\mathring{q}$.
There exists a $\mathbf{k}$-linear map $\sigma : (\hat{A}, \hat{q})\to W$ such that
\begin{equation}\label{eq3.5}
\left.\begin{array}{c}
\sigma(\hat{1})=\mathring{1}, \qquad \sigma(\hat{q})=\mathring{q}, \qquad
\sigma(\hat{q}\hat{x}_{j_1}\cdots\hat{x}_{j_m})=\mathring{q}\mathring{x}_{j_1}\cdots\mathring{x}_{j_m},\\
\sigma(\hat{x}_{i_1}\cdots \hat{x}_{i_t}\hat{x}_{j_0}\hat{q}\hat{x}_{j_1}\cdots\hat{x}_{j_m})=
\mathring{x}_{i_1}\cdots \mathring{x}_{i_t}\mathring{x}_{j_0}\mathring{q}\mathring{x}_{j_1}\cdots\mathring{x}_{j_m},\\
\sigma(\hat{x}_{j_1}\cdots\hat{x}_{j_m})=\mathring{x}_{j_1}\cdots\mathring{x}_{j_m}
\end{array}\right\}
\end{equation}
and
\begin{equation}\label{eq3.6}
\sigma(H_{u,v})=0 \quad\mbox{for all $(u, v)^{[6-th]}$-generators $H_{u,v}$,}
\end{equation}
where $ x_{i_1}, \cdots , x_{i_t}, x_{j_0}, x_{j_1}, \cdots , x_{j_m}\in X$, $i_t\geq\cdots \geq i_1$, $j_m\geq\cdots\geq j_1\geq j_0$ and 
$t$, $m\in\mathcal{Z}_{\ge 0}$.
\end{proposition}

\bigskip
For convenience, we now define nice maps which will be used to split the proof of Proposition~\ref{pr3.2} into three steps. Let $N$ be a subspace of 
$(\hat{A}, \hat{q})$. A $\mathbf{k}$-linear map $\sigma: N\to W$ is said to be {\bf nice} if $\sigma$ satisfies (\ref{eq3.5}) and (\ref{eq3.6}) for all model monomials in $N$ and for all $(u, v)^{[6-th]}$-generators in $N$. 

Recall that
$$
(\hat{A}, \hat{q})=\displaystyle\bigoplus_{m\in\mathcal{Z}_{\ge 0}}\hat{A}_m,\quad
\mbox{and}\quad \hat{A}_m=\displaystyle\bigoplus_{t=0}^{m+1}\hat{A}_m^{(t, m-t)},
$$
where
\begin{eqnarray*}
\hat{A}_m:&=&\mbox{\it the subspace spanned by all monomials whose}\\
&&\mbox{\it $X$-degrees are $m$},
\end{eqnarray*}
and
\begin{eqnarray*}
\hat{A}_m^{(t, m-t)}:&=&\mbox{\it the subspace spanned by all monomials whose}\\
&&\mbox{\it $X$-degrees are $m$ and whose $\hat{q}$-degrees are $(t, m-t)$}.
\end{eqnarray*}

Since the set
\begin{equation}\label{eq3.7}
\{\, \hat{1}, \; \hat{q}, \; \hat{q}\hat{x}_j, \; \hat{x}_j\hat{q}, \; \hat{x}_j \,|\, x_j\in X \,\}
\end{equation}
is a basis of $\hat{A}_0\oplus \hat{A}_1$, we can define a $\mathbf{k}$-linear map
$\sigma : \hat{A}_0\oplus \hat{A}_1\to W$ by
$$
\sigma (\hat{1})=\mathring{1}, \; \sigma (\hat{q})=\mathring{q}, \;
\sigma (\hat{q}\hat{x}_j)=\mathring{q}\mathring{x}_j, \; 
\sigma (\hat{x}_j\hat{q})=\mathring{x}_j\mathring{q}, \; \sigma (\hat{x}_j)=\mathring{x}_j
$$
for $x_j\in X$. The $\mathbf{k}$-linear map 
$\sigma : \hat{A}_0\oplus \hat{A}_1\to W$ is clearly nice because 
$\hat{A}_0\oplus \hat{A}_1$ does not contain any $(u, v)^{[6-th]}$-generator and the set (\ref{eq3.7}) is the set of all model monomials in $\hat{A}_0\oplus \hat{A}_1$. This proves that
\begin{eqnarray}\label{eq3.8}
&&\mbox{\it there exists a nice $\mathbf{k}$-linear map $\sigma :
\hat{A}_0\oplus \hat{A}_1\to W$.}
\end{eqnarray}

We are going to use the remaining three subsections of this section to extend the nice $\mathbf{k}$-linear map in 
(\ref{eq3.8}) to a nice  $\mathbf{k}$-linear map 
$\sigma : (\hat{A}, \hat{q})\to W$. 

\medskip
\subsection{Step 1}

The goal of Step 1 is to prove the following fact.
\begin{center}
{\it  {\bf Fact 1:} For $m\in\mathcal{Z}_{\ge 2}$, a nice $\mathbf{k}$-linear map}\\
\begin{equation}\label{eq3.9}
\sigma : \displaystyle\bigoplus_{i=0}^{m-1}\hat{A}_i\to W
\end{equation}\\
{\it can be extended to a nice $\mathbf{k}$-linear map}\\
\begin{equation}\label{eq3.10}
\sigma : \left(\displaystyle\bigoplus_{i=0}^{m-1}\hat{A}_i\right)\bigoplus
\left(\displaystyle\bigoplus_{i=0}^{2}\hat{A}_m^{(i, m-i)}
\right)\to W.
\end{equation}
\end{center}

We begin the proof of {\bf Fact 1} by introducing the central index of a monomial whose $X$-degree is at least $2$. If $m\in\mathcal{Z}_{\ge 2}$, we define the {\bf central index} $ind_c(\hat{x}_{j_1}\cdots\hat{x}_{j_t}\hat{q}
\hat{x}_{j_{t+1}}\cdots\hat{x}_{j_m})$ by
$$
ind_c(\hat{x}_{j_1}\cdots\hat{x}_{j_t}\hat{q}
\hat{x}_{j_{t+1}}\cdots
\hat{x}_{j_m}):=\left\{
\begin{array}{ll}\displaystyle\sum_{m\geq k>i\geq 1}\eta_{ik},&\quad\mbox{if $t=0$}\\&\\
\displaystyle\sum_{m\geq k>i\geq t}\eta_{ik}&\quad\mbox{if $t\ge 1$}.\end{array}\right.
$$

Let $\hat{B}^{<3}_{m,n}$ be the subspace spanned by all monomials whose $X$-degrees are $m$ with $m\in\mathcal{Z}_{\ge 2}$, whose $\hat{q}$-degrees are 
$(t, m-t)$ with $2\geq t\geq 0$, and whose central indices are $n$. Then we have
\begin{equation}\label{eq3.11}
\displaystyle\bigoplus_{i=0}^{2}\hat{A}_m^{(i, m-i)}
=\displaystyle\bigoplus_{n\in\mathcal{Z}_{\ge 0}}\hat{B}^{<3}_{m,n}.
\end{equation}
Note that every monomial in $\hat{B}^{<3}_{m,0}$ is a model monomial. Hence, if we define 
$$
\sigma\left(\hat{x}_{j_\epsilon}^{\epsilon}
\hat{x}_{j_{\epsilon+\lambda}}^{\lambda}\hat{q}
\hat{x}_{j_{\epsilon+\lambda+1}}\cdots\hat{x}_{j_m}\right):=
\mathring{x}_{j_\epsilon}^{\epsilon}
\mathring{x}_{j_{\epsilon+\lambda}}^{\lambda}\mathring{q}
\mathring{x}_{j_{\epsilon+\lambda+1}}\cdots\mathring{x}_{j_m}
$$
for $\hat{x}_{j_\epsilon}^{\epsilon}
\hat{x}_{j_{\epsilon+\lambda}}^{\lambda}\hat{q}
\hat{x}_{j_{\epsilon+\lambda+1}}\cdots\hat{x}_{j_m}\in \hat{B}^{<3}_{m,0}$ with $\epsilon$, $\lambda=0$, $1$, then the nice $\mathbf{k}$-linear map in 
(\ref{eq3.9}) is extended to a $\mathbf{k}$-linear map
\begin{equation}\label{eq3.12}
\sigma : \left(\displaystyle\bigoplus_{i=0}^{m-1}\hat{A}_i\right)\bigoplus
\hat{B}^{<3}_{m,0}\to W,
\end{equation}
and 
the  $\mathbf{k}$-linear map in (\ref{eq3.12}) satisfies (\ref{eq3.5}) for all model monomials in 
the subspace $\left(\displaystyle\bigoplus_{i=0}^{m-1}\hat{A}_i\right)\bigoplus\hat{B}^{<3}_{m,0}$.
Since every $(u, v)^{[6-th]}$-generator in $\left(\displaystyle\bigoplus_{i=0}^{m-1}\hat{A}_i\right)\bigoplus
\hat{B}^{<3}_{m,0}$ is in $\displaystyle\bigoplus_{i=0}^{m-1}\hat{A}_i$, the  $\mathbf{k}$-linear map in (\ref{eq3.12}) is nice. In other words, we have
\begin{eqnarray}\label{eq3.15}
&&\mbox{\it
for $m\in\mathcal{Z}_{\ge 2}$,  a nice  $\mathbf{k}$-linear map 
$\sigma : \displaystyle\bigoplus_{i=0}^{m-1}\hat{A}_i\to W$ can be extended}\nonumber\\
&&\mbox{\it to a nice $\mathbf{k}$-linear map $\sigma :\left(\displaystyle\bigoplus_{i=0}^{m-1}\hat{A}_i\right)\bigoplus
\hat{B}^{<3}_{m,0}\to W$. }
\end{eqnarray} 

We now prove that
\begin{eqnarray*}
&&\mbox{\it {\bf Fact 2}: For $m\in\mathcal{Z}_{\ge 2}$ and $n\in\mathcal{Z}_{\ge 1}$ ,  a nice  $\mathbf{k}$-linear map}
\end{eqnarray*}
\begin{equation}\label{eq3.16}
\sigma : \left(\displaystyle\bigoplus_{i=0}^{m-1}\hat{A}_i\right)\bigoplus
\left(\displaystyle\bigoplus_{i=0}^{n-1}\hat{B}^{<3}_{m,i}\right)\to W
\end{equation}
\begin{eqnarray*}
&&\mbox{\it can be extended to a nice  $\mathbf{k}$-linear map}
\end{eqnarray*}
\begin{equation}\label{eq3.17}
\sigma : \left(\displaystyle\bigoplus_{i=0}^{m-1}\hat{A}_i\right)\bigoplus
\left(\displaystyle\bigoplus_{i=0}^{n}\hat{B}^{<3}_{m,i}\right)\to W.
\end{equation}

Note that the set $S_1\cup S_2$ is a basis of the vector space 
$\hat{B}^{<3}_{m,n}$, where
$$
S_1:=\left\{\hat{q}
\hat{x}_{j_1}\cdots\hat{x}_{j_m}\,\left|\, \begin{array}{c}
x_{j_1}, \dots , x_{j_m}\in X,\\
ind_c(\hat{q}\hat{x}_{j_1}\cdots\hat{x}_{j_m})=n\end{array}\right.\right\}
$$
and
$$
S_2:=\left\{\hat{x}_{j_\epsilon}^{\epsilon}
\hat{x}_{j_{\epsilon+1}}\hat{q}
\hat{x}_{j_{\epsilon+2}}\cdots\hat{x}_{j_m}\left|\, \begin{array}{c}
x_{j_{\epsilon}}, x_{j_{\epsilon+1}}, \dots , x_{j_m}\in X,\\
ind_c(\hat{x}_{j_\epsilon}^{\epsilon}
\hat{x}_{j_{\epsilon+1}}\hat{q}
\hat{x}_{j_{\epsilon+2}}\cdots\hat{x}_{j_m})=n,\\
\epsilon =0, \, 1
\end{array}\right.\right\}.
$$

Using the basis $S_1\cup S_2$ of the vector space 
$\hat{B}^{<3}_{m,n}$, we define a $\mathbf{k}$-linear map 
$\sigma : \hat{B}^{<3}_{m,n}\to W$ as follows. For 
$\hat{q}\hat{x}_{j_1}\cdots\hat{x}_{j_m}\in S_1$, there exists an integer $s$ such that 
$m\geq s+1>s\geq 1$ and $j_s>j_{s+1}$, in which case, we define
\begin{eqnarray}\label{eq3.18}
&&\sigma\Big(\hat{q}\hat{x}_{j_1}\cdots
\hat{x}_{j_{s-1}}\cdot\hat{x}_{s}\hat{x}_{j_{s+1}}\cdot\hat{x}_{j_{s+2}}
\cdots\hat{x}_{j_m}\Big):=\nonumber\\
&=&\frac{1}{k}\sigma\Big(\underbrace{\hat{q}\hat{x}_{j_1}\cdots
\hat{x}_{j_{s-1}}\cdot[x_{s}, x_{j_{s+1}}]\,\hat{}\cdot\hat{x}_{j_{s+2}}
\cdots\hat{x}_{j_m}}_{(\ref{eq3.18})_1}\Big)+\nonumber\\
&&+\sigma\Big(\underbrace{\hat{q}\hat{x}_{j_1}\cdots
\hat{x}_{j_{s-1}}\cdot\hat{x}_{j_{s+1}}\hat{x}_{s}\cdot\hat{x}_{j_{s+2}}
\cdots\hat{x}_{j_m}}_{(\ref{eq3.18})_2}\Big).
\end{eqnarray}
For $\hat{x}_{j_\epsilon}^{\epsilon}
\hat{x}_{j_{\epsilon+1}}\hat{q}
\hat{x}_{j_{\epsilon+2}}\cdots\hat{x}_{j_m}\in S_2$ with $\epsilon =0, \, 1$, we have either $j_{\epsilon+1}>j_{\epsilon+2}$ or 
$j_{\epsilon+2}\geq j_{\epsilon+1}$. If $j_{\epsilon+1}>j_{\epsilon+2}$, we define
\begin{eqnarray}\label{eq3.19}
&&\sigma\Big(\hat{x}_{j_\epsilon}^{\epsilon}\cdot
\hat{x}_{j_{\epsilon+1}}\hat{q}
\hat{x}_{j_{\epsilon+2}}\cdot\hat{x}_{j_{\epsilon+3}}\cdots\hat{x}_{j_m}\Big)
:=\nonumber\\
&=&\frac{1}{k}\sigma\Big(\underbrace{\hat{x}_{j_\epsilon}^{\epsilon}\cdot
[x_{j_{\epsilon+1}},
x_{j_{\epsilon+2}}]\,\hat{}\cdot\hat{q}\hat{x}_{j_{\epsilon+3}}
\cdots\hat{x}_{j_m}}_{(\ref{eq3.19})_1}\Big)+\nonumber\\
&&+\sigma\Big(\underbrace{\hat{x}_{j_\epsilon}^{\epsilon}\cdot
\hat{x}_{j_{\epsilon+2}}\hat{q}
\hat{x}_{j_{\epsilon+1}}\cdot\hat{x}_{j_{\epsilon+3}}
\cdots\hat{x}_{j_m}}_{(\ref{eq3.19})_2}\Big);
\end{eqnarray}
if $j_{\epsilon+2}\geq j_{\epsilon+1}$, there exists an integer $s$ such that
$m\geq s+1>s\geq \epsilon+2$ and $j_s>j_{s+1}$, in which case, we define
\begin{eqnarray}\label{eq3.20}
&&\sigma\Big(\hat{x}_{j_\epsilon}^{\epsilon}
\hat{x}_{j_{\epsilon+1}}\hat{q}
\hat{x}_{j_{\epsilon+2}}\cdots
\hat{x}_{j_{s-1}}\cdot\hat{x}_{j_s}\hat{x}_{j_{s+1}}
\cdot\hat{x}_{j_{s+2}}\cdots\hat{x}_{j_m}\Big):=\nonumber\\
&=&\frac{1}{k}\sigma\Big(\underbrace{\hat{x}_{j_\epsilon}^{\epsilon}
\hat{x}_{j_{\epsilon+1}}\hat{q}
\hat{x}_{j_{\epsilon+2}}\cdots
\hat{x}_{j_{s-1}}\cdot[x_{j_s}, x_{j_{s+1}}]\,\hat{}
\cdot\hat{x}_{j_{s+2}}\cdots\hat{x}_{j_m}}_{(\ref{eq3.20})_1}\Big)+\nonumber\\
&&+\sigma\Big(\underbrace{\hat{x}_{j_\epsilon}^{\epsilon}
\hat{x}_{j_{\epsilon+1}}\hat{q}
\hat{x}_{j_{\epsilon+2}}\cdots
\hat{x}_{j_{s-1}}\cdot\hat{x}_{j_{s+1}}\hat{x}_{j_s}
\cdot\hat{x}_{j_{s+2}}\cdots\hat{x}_{j_m}}_{(\ref{eq3.20})_2}\Big).
\end{eqnarray}
Since 
$$
deg_X((\ref{eq3.18})_1)=deg_X((\ref{eq3.19})_1)
=deg_X((\ref{eq3.20})_1)=m-1<m
$$
and
$$
ind_c((\ref{eq3.18})_2)=ind_c((\ref{eq3.19})_2)=ind_c((\ref{eq3.20})_2)=n-1<n,
$$
each of the right hand sides of (\ref{eq3.18}), (\ref{eq3.19}) and (\ref{eq3.20}) makes sense by using the $\mathbf{k}$-linear map $\sigma$ in (\ref{eq3.16}).

Although (\ref{eq3.19}) is well-defined, we need to prove that both (\ref{eq3.18}) and (\ref{eq3.20}) are well defined. We shall just prove that 
(\ref{eq3.20}) is well defined because a similar argument can be used to prove that (\ref{eq3.18}) is also well defined. 

Suppose that there exists another pair $(j_t, j_{t+1})$ in (\ref{eq3.20}) such that 
$m\geq t+1>t\geq \epsilon +2$ and $j_t>j_{t+1}$. Using the pair 
$(j_t, j_{t+1})$ and the definition given by (\ref{eq3.20}), we have
\begin{eqnarray}\label{eq3.21}
&&\sigma\Big(\hat{x}_{j_\epsilon}^{\epsilon}
\hat{x}_{j_{\epsilon+1}}\hat{q}
\hat{x}_{j_{\epsilon+2}}\cdots
\hat{x}_{j_{t-1}}\cdot\hat{x}_{j_t}\hat{x}_{j_{t+1}}
\cdot\hat{x}_{j_{t+2}}\cdots\hat{x}_{j_m}\Big):=\nonumber\\
&=&\frac{1}{k}\sigma\Big(\underbrace{\hat{x}_{j_\epsilon}^{\epsilon}
\hat{x}_{j_{\epsilon+1}}\hat{q}
\hat{x}_{j_{\epsilon+2}}\cdots
\hat{x}_{j_{t-1}}\cdot[x_{j_t}, x_{j_{t+1}}]\,\hat{}
\cdot\hat{x}_{j_{t+2}}\cdots\hat{x}_{j_m}}_{(\ref{eq3.21})_1}\Big)+\nonumber\\
&&+\sigma\Big(\underbrace{\hat{x}_{j_\epsilon}^{\epsilon}
\hat{x}_{j_{\epsilon+1}}\hat{q}
\hat{x}_{j_{\epsilon+2}}\cdots
\hat{x}_{j_{t-1}}\cdot\hat{x}_{j_{t+1}}\hat{x}_{j_t}
\cdot\hat{x}_{j_{t+2}}\cdots\hat{x}_{j_m}}_{(\ref{eq3.21})_2}\Big).
\end{eqnarray}

\medskip
One can prove that
\begin{equation}\label{eq3.22}
\mbox{the right hand side of (\ref{eq3.20})}=
\mbox{the right hand side of (\ref{eq3.21})}.
\end{equation}

\medskip

Using (\ref{eq3.18}), (\ref{eq3.19}) and (\ref{eq3.20}), the $\mathbf{k}$-linear map 
in (\ref{eq3.16}) can be extended to a  $\mathbf{k}$-linear map in (\ref{eq3.17}). To finish Step 1, we need only to prove that
\begin{eqnarray}\label{eq3.37}
&&\mbox{\it
If  $\sigma :  \left(\displaystyle\bigoplus_{i=0}^{m-1}\hat{A}_i\right)\bigoplus
\left(\displaystyle\bigoplus_{i=0}^{n}\hat{B}^{<3}_{m,i}\right)\to W$ is a  $\mathbf{k}$-linear map which}\nonumber\\
&&\mbox{\it extends the nice $\mathbf{k}$-linear map 
in (\ref{eq3.16}) and satisfies (\ref{eq3.18}),}\nonumber\\
&&\mbox{\it (\ref{eq3.19}) and (\ref{eq3.20}), then $\sigma$ is nice.}
\end{eqnarray} 

Note that any $(u, v)^{[6-th]}$-generator $H^{<3}_{uv}$ satisfying
\begin{equation}\label{eq3.38}
H^{<3}_{uv}\in \left(\displaystyle\bigoplus_{i=0}^{m-1}\hat{A}_i\right)\bigoplus
\left(\displaystyle\bigoplus_{i=0}^{n}\hat{B}^{<3}_{m,i}\right)
\end{equation}
and
\begin{equation}\label{eq3.39}
H^{<3}_{uv}\not\in \left(\displaystyle\bigoplus_{i=0}^{m-1}\hat{A}_i\right)\bigoplus
\left(\displaystyle\bigoplus_{i=0}^{n-1}\hat{B}^{<3}_{m,i}\right)
\end{equation}
is either a $(1, 0)^{[6-th]}$-generator or a $(0, 1)^{[6-th]}$-generator. 

Let $H^{<3}_{10}$ be a $(1, 0)^{[6-th]}$-generator satisfying (\ref{eq3.38}) and (\ref{eq3.39}). Then we have
\begin{eqnarray*}
H^{<3}_{10}&=&\underbrace{\hat{x}_{j_\epsilon}^{\epsilon}
\hat{x}_{j_{\epsilon+1}}\hat{q}
\hat{x}_{j_{\epsilon+2}}\hat{x}_{j_{\epsilon+3}}\cdots\hat{x}_{j_{s-1}}\cdot
[x_{j_s}, x_{j_{s+1}}]\,\hat{}\cdot\hat{x}_{j_{s+2}}\cdots\hat{x}_{j_m}}_{H^{<3}_{10, 1}}+\\
&&-k\underbrace{\hat{x}_{j_\epsilon}^{\epsilon}
\hat{x}_{j_{\epsilon+1}}\hat{q}
\hat{x}_{j_{\epsilon+2}}\hat{x}_{j_{\epsilon+3}}\cdots\hat{x}_{j_{s-1}}
(\hat{x}_{j_s}\hat{x}_{j_{s+1}})\hat{x}_{j_{s+2}}\cdots\hat{x}_{j_m}}_{H^{<3}_{10, 2}}+\\
&&+k\underbrace{\hat{x}_{j_\epsilon}^{\epsilon}
\hat{x}_{j_{\epsilon+1}}\hat{q}
\hat{x}_{j_{\epsilon+2}}\hat{x}_{j_{\epsilon+3}}\cdots\hat{x}_{j_{s-1}}
(\hat{x}_{j_{s+1}}\hat{x}_{j_s})\hat{x}_{j_{s+2}}\cdots\hat{x}_{j_m}}_{H^{<3}_{10, 3}}
\end{eqnarray*}
by (\ref{eq3.3}), where $j_s\ne j_{s+1}$ and $\hat{x}_{j_\epsilon}^{\epsilon}
\hat{x}_{j_{\epsilon+1}}$ maybe disappear. Without losing the generality, we can assume that 
$j_s> j_{s+1}$. Hence, $ind_c\Big(H^{<3}_{10, 2}\Big)=n$. If $j_{\epsilon+2}\geq j_{\epsilon+1}$ or 
$\hat{x}_{j_\epsilon}^{\epsilon}\hat{x}_{j_{\epsilon+1}}$ disappears, then 
$\sigma\Big(H^{<3}_{10}\Big)=0$ by (\ref{eq3.20}) or (\ref{eq3.18}). If 
$j_{\epsilon+1}>j_{\epsilon+2}$, then (\ref{eq3.19}) can be used to rewrite $H^{<3}_{10}$ as follows:
\begin{eqnarray}
\sigma\Big(H^{<3}_{10}\Big)&=&\sigma\Big(\underbrace{\hat{x}_{j_\epsilon}^{\epsilon}
(\hat{x}_{j_{\epsilon+1}}\hat{q}
\hat{x}_{j_{\epsilon+2}})\hat{x}_{j_{\epsilon+3}}\cdots\hat{x}_{j_{s-1}}\cdot
[x_{j_s}, x_{j_{s+1}}]\,\hat{}\cdot\hat{x}_{j_{s+2}}\cdots\hat{x}_{j_m}}_{H^{<3}_{10, 1}}\Big)+\nonumber\\
&&-\sigma\Big(\hat{x}_{j_\epsilon}^{\epsilon}
[x_{j_{\epsilon+1}}, x_{j_{\epsilon+2}}]\,\hat{}\;\hat{q}
\hat{x}_{j_{\epsilon+3}}\cdots\hat{x}_{j_{s-1}}
(\hat{x}_{j_s}\hat{x}_{j_{s+1}})\hat{x}_{j_{s+2}}\cdots\hat{x}_{j_m}\Big)+\nonumber\\
&&-k\sigma\Big(\hat{x}_{j_\epsilon}^{\epsilon}
x_{j_{\epsilon+2}}\hat{q}x_{j_{\epsilon+1}}
\hat{x}_{j_{\epsilon+3}}\cdots\hat{x}_{j_{s-1}}
(\hat{x}_{j_s}\hat{x}_{j_{s+1}})\hat{x}_{j_{s+2}}\cdots\hat{x}_{j_m}\Big)+\nonumber
\end{eqnarray} 
\begin{equation}\label{eq3.40}
+k\sigma\Big(\underbrace{\hat{x}_{j_\epsilon}^{\epsilon}
(\hat{x}_{j_{\epsilon+1}}\hat{q}
\hat{x}_{j_{\epsilon+2}})\hat{x}_{j_{\epsilon+3}}\cdots\hat{x}_{j_{s-1}}
(\hat{x}_{j_{s+1}}\hat{x}_{j_s})\hat{x}_{j_{s+2}}\cdots\hat{x}_{j_m}}_{H^{<3}_{10, 3}}\Big).
\end{equation}

Let
\begin{eqnarray}
H_{01}'&=&\hat{x}_{j_\epsilon}^{\epsilon}\cdot
[x_{j_{\epsilon+1}}, x_{j_{\epsilon+2}}]\,\hat{}\;\cdot\hat{q}\hat{x}_{j_{\epsilon+3}}\cdots
\hat{x}_{j_{s-1}}\cdot[x_{j_s}, x_{j_{s+1}}]\,\hat{}
\cdot\hat{x}_{j_{s+2}}\cdots\hat{x}_{j_m}+\nonumber\\
&&-k\underbrace{\hat{x}_{j_\epsilon}^{\epsilon}
(\hat{x}_{j_{\epsilon+1}}\hat{q}
\hat{x}_{j_{\epsilon+2}})\hat{x}_{j_{\epsilon+3}}\cdots\hat{x}_{j_{s-1}}\cdot
[x_{j_s}, x_{j_{s+1}}]\,\hat{}\;\cdot\hat{x}_{j_{s+2}}
\cdots\hat{x}_{j_m}}_{H^{<3}_{10, 1}}+\nonumber
\end{eqnarray}
\begin{equation}\label{eq3.41}
+k\hat{x}_{j_\epsilon}^{\epsilon}
(\hat{x}_{j_{\epsilon+2}}\hat{q}
\hat{x}_{j_{\epsilon+1}})\hat{x}_{j_{\epsilon+3}}\cdots\hat{x}_{j_{s-1}}\cdot
[x_{j_s}, x_{j_{s+1}}]\,\hat{}\;\cdot\hat{x}_{j_{s+2}}
\cdots\hat{x}_{j_m}
\end{equation} 
and
\begin{eqnarray}
H_{01}''&=&\hat{x}_{j_\epsilon}^{\epsilon}\cdot
[x_{j_{\epsilon+1}}, x_{j_{\epsilon+2}}]\,\hat{}\;\cdot\hat{q}\hat{x}_{j_{\epsilon+3}}\cdots
\hat{x}_{j_{s-1}}(\hat{x}_{j_{s+1}}\hat{x}_{j_s})
\hat{x}_{j_{s+2}}\cdots\hat{x}_{j_m}+\nonumber\\
&&-k\underbrace{\hat{x}_{j_\epsilon}^{\epsilon}
(\hat{x}_{j_{\epsilon+1}}\hat{q}
\hat{x}_{j_{\epsilon+2}})\hat{x}_{j_{\epsilon+3}}\cdots\hat{x}_{j_{s-1}}
(\hat{x}_{j_{s+1}}\hat{x}_{j_s})
\hat{x}_{j_{s+2}}
\cdots\hat{x}_{j_m}}_{H^{<3}_{10, 3}}+\nonumber
\end{eqnarray}
\begin{equation}\label{eq3.42}
+k\hat{x}_{j_\epsilon}^{\epsilon}
(\hat{x}_{j_{\epsilon+2}}\hat{q}
\hat{x}_{j_{\epsilon+1}})\hat{x}_{j_{\epsilon+3}}\cdots\hat{x}_{j_{s-1}}
(\hat{x}_{j_{s+1}}\hat{x}_{j_s})\hat{x}_{j_{s+2}}
\cdots\hat{x}_{j_m}
\end{equation} 
Then $H_{01}'$ is a linear combination of 
$(0, 1)^{[6-th]}$-generators and $H_{01}''$ is a $(0, 1)^{[6-th]}$-generator. Clearly, $H_{01}'\in \displaystyle\bigoplus_{i=0}^{m-1}\hat{A}_i$. Since $j_s>j_{s+1}$ and
$ind_c\Big(H^{<3}_{10, 2}\Big)=n$, $ind_c\Big(H^{<3}_{10, 3}\Big)=n-1$. Thus, 
$H_{01}''\in \left(\displaystyle\bigoplus_{i=0}^{m-1}\hat{A}_i\right)\bigoplus
\left(\displaystyle\bigoplus_{i=0}^{n-1}\hat{B}^{<3}_{m,i}\right)$. Using the nice 
$\mathbf{k}$-linear map in (\ref{eq3.16}), we get
\begin{equation}\label{eq3.43}
\sigma\Big(H_{01}'\Big)=0\quad\mbox{and}\quad \sigma\Big(H_{01}''\Big)=0.
\end{equation}

It follows from (\ref{eq3.41}), (\ref{eq3.42}) and (\ref{eq3.43}) that
\begin{eqnarray*}
&&\sigma\Big(\underbrace{\hat{x}_{j_\epsilon}^{\epsilon}
(\hat{x}_{j_{\epsilon+1}}\hat{q} x_{j_{\epsilon+2}})\hat{x}_{j_{\epsilon+3}}\cdots
\hat{x}_{j_{s-1}}\cdot[x_{j_s}, x_{j_{s+1}}]\,\hat{}
\cdot\hat{x}_{j_{s+2}}\cdots\hat{x}_{j_m}}_{H^{<3}_{10, 1}}\Big)\\
&=&\frac{1}{k}\sigma\Big(\hat{x}_{j_\epsilon}^{\epsilon}
[x_{j_{\epsilon+1}}, x_{j_{\epsilon+2}}]\,\hat{}\;\hat{q}\hat{x}_{j_{\epsilon+3}}\cdots
\hat{x}_{j_{s-1}}\cdot[x_{j_s}, x_{j_{s+1}}]\,\hat{}
\cdot\hat{x}_{j_{s+2}}\cdots\hat{x}_{j_m}\Big)+\\
\end{eqnarray*}
\begin{equation}\label{eq3.44}
+\sigma\Big(\hat{x}_{j_\epsilon}^{\epsilon}
(\hat{x}_{j_{\epsilon+2}}\hat{q} x_{j_{\epsilon+1}})\hat{x}_{j_{\epsilon+3}}\cdots
\hat{x}_{j_{s-1}}\cdot[x_{j_s}, x_{j_{s+1}}]\,\hat{}
\cdot\hat{x}_{j_{s+2}}\cdots\hat{x}_{j_m}\Big)
\end{equation} 
and
\begin{eqnarray*}
&&\sigma\Big(\underbrace{\hat{x}_{j_\epsilon}^{\epsilon}
(\hat{x}_{j_{\epsilon+1}}\hat{q} x_{j_{\epsilon+2}})\hat{x}_{j_{\epsilon+3}}\cdots
\hat{x}_{j_{s-1}}(\hat{x}_{j_{s+1}}\hat{x}_{j_s})
\hat{x}_{j_{s+2}}\cdots\hat{x}_{j_m}}_{H^{<3}_{10, 3}}\Big)\\
&=&\frac{1}{k}\sigma\Big(\hat{x}_{j_\epsilon}^{\epsilon}
[x_{j_{\epsilon+1}}, x_{j_{\epsilon+2}}]\,\hat{}\;\hat{q}\hat{x}_{j_{\epsilon+3}}\cdots
\hat{x}_{j_{s-1}}(\hat{x}_{j_{s+1}}\hat{x}_{j_s})
\hat{x}_{j_{s+2}}\cdots\hat{x}_{j_m}\Big)+\\
\end{eqnarray*}
\begin{equation}\label{eq3.45}
+\sigma\Big(\hat{x}_{j_\epsilon}^{\epsilon}
(\hat{x}_{j_{\epsilon+2}}\hat{q} x_{j_{\epsilon+1}})\hat{x}_{j_{\epsilon+3}}\cdots
\hat{x}_{j_{s-1}}(\hat{x}_{j_{s+1}}\hat{x}_{j_s})
\hat{x}_{j_{s+2}}\cdots\hat{x}_{j_m}\Big).
\end{equation} 

By (\ref{eq3.44}) and (\ref{eq3.45}), (\ref{eq3.40}) becomes
\begin{eqnarray*}
&&\sigma\Big(H^{<3}_{10}\Big)\\
&=&\underbrace{\frac{1}{k}\sigma\Big(\hat{x}_{j_\epsilon}^{\epsilon}
[x_{j_{\epsilon+1}}, x_{j_{\epsilon+2}}]\,\hat{}\;\hat{q}\hat{x}_{j_{\epsilon+3}}\cdots
\hat{x}_{j_{s-1}}\cdot[x_{j_s}, x_{j_{s+1}}]\,\hat{}
\cdot\hat{x}_{j_{s+2}}\cdots\hat{x}_{j_m}\Big)}_{(\ref{eq3.46})_1}+\\
&&+\underbrace{\sigma\Big(\hat{x}_{j_\epsilon}^{\epsilon}
(\hat{x}_{j_{\epsilon+2}}\hat{q} x_{j_{\epsilon+1}})\hat{x}_{j_{\epsilon+3}}\cdots
\hat{x}_{j_{s-1}}\cdot[x_{j_s}, x_{j_{s+1}}]\,\hat{}
\cdot\hat{x}_{j_{s+2}}\cdots\hat{x}_{j_m}\Big)}_{(\ref{eq3.46})_2}\\
&&-\underbrace{\sigma\Big(\hat{x}_{j_\epsilon}^{\epsilon}
[x_{j_{\epsilon+1}}, x_{j_{\epsilon+2}}]\,\hat{}\;\hat{q}
\hat{x}_{j_{\epsilon+3}}\cdots\hat{x}_{j_{s-1}}
(\hat{x}_{j_s}\hat{x}_{j_{s+1}})\hat{x}_{j_{s+2}}\cdots\hat{x}_{j_m}\Big)}_{(\ref{eq3.46})_3}+\\
&&-\underbrace{k\sigma\Big(\hat{x}_{j_\epsilon}^{\epsilon}
x_{j_{\epsilon+2}}\hat{q}x_{j_{\epsilon+1}}
\hat{x}_{j_{\epsilon+3}}\cdots\hat{x}_{j_{s-1}}
(\hat{x}_{j_s}\hat{x}_{j_{s+1}})\hat{x}_{j_{s+2}}\cdots\hat{x}_{j_m}\Big)}_{(\ref{eq3.46})_4}+\\
&&+\underbrace{\sigma\Big(\hat{x}_{j_\epsilon}^{\epsilon}
[x_{j_{\epsilon+1}}, x_{j_{\epsilon+2}}]\,\hat{}\;\hat{q}\hat{x}_{j_{\epsilon+3}}\cdots
\hat{x}_{j_{s-1}}(\hat{x}_{j_{s+1}}\hat{x}_{j_s})
\hat{x}_{j_{s+2}}\cdots\hat{x}_{j_m}\Big)}_{(\ref{eq3.46})_5}+
\end{eqnarray*}
\begin{equation}\label{eq3.46}
+\underbrace{k\sigma\Big(\hat{x}_{j_\epsilon}^{\epsilon}
(\hat{x}_{j_{\epsilon+2}}\hat{q} x_{j_{\epsilon+1}})\hat{x}_{j_{\epsilon+3}}\cdots
\hat{x}_{j_{s-1}}(\hat{x}_{j_{s+1}}\hat{x}_{j_s})
\hat{x}_{j_{s+2}}\cdots\hat{x}_{j_m}\Big)}_{(\ref{eq3.46})_6}.
\end{equation}

Note that
\begin{eqnarray}\label{eq3.47}
&&(\ref{eq3.46})_1-(\ref{eq3.46})_3+(\ref{eq3.46})_5\nonumber\\
&=&k\sigma\left(\mbox{\it a $(1, 0)^{[6-th]}$-generator in 
$\displaystyle\bigoplus_{i=0}^{m-1}\hat{A}_i$}\right)=0
\end{eqnarray}
and
\begin{eqnarray}\label{eq3.48}
&&(\ref{eq3.46})_2-(\ref{eq3.46})_4+(\ref{eq3.46})_6\\
&=&\sigma\left(\mbox{\it a $(1, 0)^{[6-th]}$-generator in 
$\left(\displaystyle\bigoplus_{i=0}^{m-1}\hat{A}_i\right)\bigoplus
\left(\displaystyle\bigoplus_{i=0}^{n-1}\hat{B}^{<3}_{m,i}\right)$}\right)=0.\nonumber
\end{eqnarray}
It follows from (\ref{eq3.46}), (\ref{eq3.47}) and (\ref{eq3.48}) that
\begin{eqnarray}\label{eq3.49}
&&\mbox{\it If $H^{<3}_{10}$ is a $(1, 0)^{[6-th]}$-generator satisfying (\ref{eq3.38}) and }\nonumber\\
&&\mbox{(\ref{eq3.39}), then $\sigma\Big(H^{<3}_{10}\Big)=0$.}
\end{eqnarray}

Let $H^{<3}_{01}$ be a $(0, 1)^{[6-th]}$-generator satisfying (\ref{eq3.38}) and (\ref{eq3.39}), then we have
\begin{eqnarray*}
H^{<3}_{01}&=&[x_{j_1}, x_{j_2}]\,\hat{}
\cdot\hat{q}\hat{x}_{j_3}\cdots\hat{x}_{j_m}+\\
&&-k(\hat{x}_{j_1}\hat{q}\hat{x}_{j_2})\hat{x}_{j_3}\cdots\hat{x}_{j_m}+
k(\hat{x}_{j_2}\hat{q}\hat{x}_{j_1})\hat{x}_{j_3}\cdots\hat{x}_{j_m}
\end{eqnarray*}
by (\ref{eq3.4}), where $j_1\ne j_2$. Clearly, $\sigma\Big(H^{<3}_{01}\Big)=0$ by (\ref{eq3.19}). This proves that
\begin{eqnarray}\label{eq3.50}
&&\mbox{\it If $H^{<3}_{01}$ is a $(0, 1)^{[6-th]}$-generator satisfying (\ref{eq3.38}) and }\nonumber\\
&&\mbox{(\ref{eq3.39}), then $\sigma\Big(H^{<3}_{01}\Big)=0$.}
\end{eqnarray}

By (\ref{eq3.49}) and (\ref{eq3.50}), (\ref{eq3.37}) holds. Using (\ref{eq3.11}), (\ref{eq3.15}) and (\ref{eq3.16}), the nice $\mathbf{k}$-linear map in (\ref{eq3.9}) can be extended a nice $\mathbf{k}$-linear map in (\ref{eq3.10}). This completes Step 1.

\medskip
\subsection{Step 2}

The goal of Step 2 is to prove the following fact.
\begin{center}
{\it  {\bf Fact 3:} For $m\in\mathcal{Z}_{\ge 3}$ and $m\geq t\geq 3$, 
a nice $\mathbf{k}$-linear map}\\
\begin{equation}\label{eq3.51}
\sigma : \left(\displaystyle\bigoplus_{i=0}^{m-1}\hat{A}_i\right)\bigoplus
\left(\displaystyle\bigoplus_{i=0}^{t-1}\hat{A}_m^{(i, m-i)}
\right)\to W
\end{equation}\\
{\it can be extended to a nice $\mathbf{k}$-linear map}\\
\begin{equation}\label{eq3.52}
\sigma : \left(\displaystyle\bigoplus_{i=0}^{m-1}\hat{A}_i\right)\bigoplus
\left(\displaystyle\bigoplus_{i=0}^{t}\hat{A}_m^{(i, m-i)}
\right)\to W.
\end{equation}
\end{center}

Let $\hat{C}^{\ge 3}_{m,t,n}$ be the subspace spanned by all monomials whose $X$-degrees are $m$ with $m\in\mathcal{Z}_{\ge 2}$, whose $\hat{q}$-degrees are 
$(t, m-t)$ with $m\ge t\ge 3$, and whose left indices are $n$. Then we have
\begin{equation}\label{eq3.53}
\hat{A}_m^{(t, m-t)}
=\displaystyle\bigoplus_{n\in\mathcal{Z}_{\ge 0}}\hat{C}^{\ge 3}_{m,t,n}.
\end{equation}

First, we prove that 
\begin{eqnarray*}
&&\mbox{\it
For $m\in\mathcal{Z}_{\ge 3}$ and $m\geq t\geq 3$,  the nice  $\mathbf{k}$-linear map 
$\sigma$}\nonumber\\
&&\mbox{\it  in (\ref{eq3.51}) can be extended to a nice $\mathbf{k}$-linear map}
\end{eqnarray*} 
\begin{equation}\label{eq3.54}
\sigma : \left(\displaystyle\bigoplus_{i=0}^{m-1}\hat{A}_i\right)\bigoplus
\left(\displaystyle\bigoplus_{i=0}^{t-1}\hat{A}_m^{(i, m-i)}\right)\bigoplus\hat{C}^{\ge 3}_{m,t,0}
\to W. 
\end{equation} 

By the definition of $\hat{C}^{\ge 3}_{m,t,0}$, we know that
$$
\left\{\, \hat{x}_{j_1}\cdots\hat{x}_{j_{t-1}}\hat{x}_{j_t}\hat{q}
\hat{x}_{j_{t+1}}\cdots\hat{x}_{j_m}\left|\, \begin{array}{c}
x_{j_1}, \dots , x_{j_m}\in X,\\
j_{t-1}\ge \cdots\ge j_1
\end{array}\right.\right\}
$$
is a basis of $\hat{C}^{\ge 3}_{m,t,0}$. Let $\hat{C}^{\ge 3}_{m,t,0, n}$ be the subspace spanned by all monomials whose $X$-degrees are $m$ with $m\in\mathcal{Z}_{\ge 3}$, whose $\hat{q}$-degrees are 
$(t, m-t)$ with $m\geq t\geq 3$, whose left indices are $0$, and whose central indices are $n$. Then we have
\begin{equation}\label{eq3.55}
\hat{C}^{\ge 3}_{m,t,0}=\bigoplus_{n\in\mathcal{Z}_{\ge 0}}\hat{C}^{\ge 3}_{m,t,0, n}.
\end{equation}

Note that every monomial in $\hat{C}^{\ge 3}_{m,t,0, 0}$ is a model monomial. Hence, if we define
$$
\sigma( \hat{x}_{j_1}\cdots\hat{x}_{j_{t-1}}\hat{x}_{j_t}\hat{q}
\hat{x}_{j_{t+1}}\cdots\hat{x}_{j_m}):= \hat{x}_{j_1}\cdots\hat{x}_{j_{t-1}}\hat{x}_{j_t}\hat{q}
\hat{x}_{j_{t+1}}\cdots\hat{x}_{j_m}
$$
for $ \hat{x}_{j_1}\cdots\hat{x}_{j_{t-1}}\hat{x}_{j_t}\hat{q}
\hat{x}_{j_{t+1}}\cdots\hat{x}_{j_m}\in \hat{C}^{\ge 3}_{m,t,0, 0}$, then the nice $\mathbf{k}$-linear map $\sigma$ in 
(\ref{eq3.51}) can be extended to a nice $\mathbf{k}$-linear map
\begin{equation}\label{eq3.56}
\sigma : \left(\displaystyle\bigoplus_{i=0}^{m-1}\hat{A}_i\right)\bigoplus
\left(\displaystyle\bigoplus_{i=0}^{t-1}\hat{A}_m^{(i, m-i)}\right)\bigoplus\hat{C}^{\ge 3}_{m,t,0, 0}\to W. 
\end{equation}
Since every  $(u, v)^{[6-th]}$-generator in $\left(\displaystyle\bigoplus_{i=0}^{m-1}\hat{A}_i\right)\bigoplus
\left(\displaystyle\bigoplus_{i=0}^{t-1}\hat{A}_m^{(i, m-i)}\right)\bigoplus\hat{C}^{\ge 3}_{m,t,0, 0}$ is in 
$\left(\displaystyle\bigoplus_{i=0}^{m-1}\hat{A}_i\right)\bigoplus
\left(\displaystyle\bigoplus_{i=0}^{t-1}\hat{A}_m^{(i, m-i)}\right)$, the  $\mathbf{k}$-linear map in 
(\ref{eq3.56}) is nice. Thus, we get
\begin{eqnarray*}
&&\mbox{\it
For $m\in\mathcal{Z}_{\ge 3}$ and $m\geq t\geq 3$,  the nice  $\mathbf{k}$-linear map 
$\sigma$}\nonumber\\
&&\mbox{\it  in (\ref{eq3.51}) can be extended to a nice $\mathbf{k}$-linear map}
\end{eqnarray*} 
\begin{equation}\label{eq3.57}
\sigma : \left(\displaystyle\bigoplus_{i=0}^{m-1}\hat{A}_i\right)\bigoplus
\left(\displaystyle\bigoplus_{i=0}^{t-1}\hat{A}_m^{(i, m-i)}\right)\bigoplus\hat{C}^{\ge 3}_{m,t,0,0}\to W. 
\end{equation} 

We now prove that
\begin{eqnarray*}
&&\mbox{\it For $m\in\mathcal{Z}_{\ge 3}$ and $n\in \mathcal{Z}_{\ge 1}$, a nice  $\mathbf{k}$-linear map }
\end{eqnarray*} 
\begin{equation}\label{eq3.58}
\sigma : \left(\displaystyle\bigoplus_{i=0}^{m-1}\hat{A}_i\right)\bigoplus
\left(\displaystyle\bigoplus_{i=0}^{t-1}\hat{A}_m^{(i, m-i)}\right)\bigoplus
\left(\displaystyle\bigoplus_{i=0}^{n-1}\hat{C}^{\ge 3}_{m,t,0,i}\right)\to W 
\end{equation} 
\begin{eqnarray*}
&&\mbox{\it can be extended to a nice  $\mathbf{k}$-linear map }
\end{eqnarray*} 
\begin{equation}\label{eq3.59}
\sigma : \left(\displaystyle\bigoplus_{i=0}^{m-1}\hat{A}_i\right)\bigoplus
\left(\displaystyle\bigoplus_{i=0}^{t-1}\hat{A}_m^{(i, m-i)}\right)\bigoplus
\left(\displaystyle\bigoplus_{i=0}^{n}\hat{C}^{\ge 3}_{m,t,0,i}\right)\to W .
\end{equation} 

Clearly, the set $S_n$ defined by
$$
\left\{ \hat{x}_{j_1}\cdots\hat{x}_{j_{t-1}}\hat{x}_{j_t}\hat{q}
\hat{x}_{j_{t+1}}\cdots\hat{x}_{j_m}\left|\, \begin{array}{c}
x_{j_1}, \dots , x_{j_m}\in X,\\
j_{t-1}\ge \cdots \ge j_1, \quad m\ge t\ge 3\\
ind_c(\hat{x}_{j_1}\cdots\hat{x}_{j_{t-1}}\hat{x}_{j_t}\hat{q}
\hat{x}_{j_{t+1}}\cdots\hat{x}_{j_m})=n
\end{array}\right.\right\}
$$
is a basis of $\hat{C}^{\ge 3}_{m,t,0, n}$. For $\hat{x}_{j_1}\cdots\hat{x}_{j_{t-1}}\hat{x}_{j_t}\hat{q}
\hat{x}_{j_{t+1}}\cdots\hat{x}_{j_m}\in S_n$, we have either $j_t>j_{t+1}$ or $j_{t+1}\ge j_t$. If $j_t>j_{t+1}$, we define
\begin{eqnarray}\label{eq3.60}
&&\sigma\Big(\hat{x}_{j_1}\cdots\hat{x}_{j_{t-1}}\hat{x}_{j_t}\hat{q}
\hat{x}_{j_{t+1}}\hat{x}_{j_{t+2}}\cdots\hat{x}_{j_m}\Big)\nonumber\\
&:=&
\frac{1}{k}\sigma\Big(\underbrace{\hat{x}_{j_1}\cdots\hat{x}_{j_{t-1}}\cdot[x_{j_t}, x_{j_{t+1}}]\,\hat{}\;\cdot
\hat{q}\hat{x}_{j_{t+2}}\cdots\hat{x}_{j_m}}_{(\ref{eq3.60})_1}\Big)+\nonumber\\
&&+\sigma\Big(\underbrace{\hat{x}_{j_1}\cdots\hat{x}_{j_{t-1}}\hat{x}_{j_{t+1}}\hat{q}
\hat{x}_{j_t}\hat{x}_{j_{t+2}}\cdots\hat{x}_{j_m}}_{(\ref{eq3.60})_2}\Big).
\end{eqnarray} 
If $j_{t+1}\ge j_t$, there exists an integer $s$ such that $m\ge s+1>s\ge t+2 $ and $j_s>j_{s+1}$, in which case, we define
\begin{eqnarray}\label{eq3.61}
&&\sigma\Big(\hat{x}_{j_1}\cdots\hat{x}_{j_{t-1}}\hat{x}_{j_t}\hat{q}
\hat{x}_{j_{t+1}}\cdots\hat{x}_{j_{s-1}}\cdot \hat{x}_{j_s}\hat{x}_{j_{s+1}}\cdot\hat{x}_{j_{s+2}}\cdots\hat{x}_{j_m}\Big)\nonumber\\
&:=&
\frac{1}{k}\sigma\Big(\underbrace{\hat{x}_{j_1}\cdots\hat{x}_{j_{t-1}}\hat{x}_{j_t}\hat{q} \hat{x}_{j_{t+1}}
\cdots\hat{x}_{j_{s-1}}\cdot[x_{j_s}, x_{j_{s+1}}]\,\hat{}\;\cdot\hat{x}_{j_{s+2}}\cdots\hat{x}_{j_m}}_{(\ref{eq3.61})_1}\Big)+\nonumber\\
&&+\sigma\Big(\underbrace{\hat{x}_{j_1}\cdots\hat{x}_{j_{t-1}}\hat{x}_{j_t}\hat{q}
\hat{x}_{j_{t+1}}\cdots\hat{x}_{j_{s-1}}\cdot \hat{x}_{j_{s+1}}\hat{x}_{j_s}\cdot\hat{x}_{j_{s+2}}
\cdots\hat{x}_{j_m}}_{(\ref{eq3.61})_2}\Big).
\end{eqnarray} 
Since $(\ref{eq3.60})_1$, $(\ref{eq3.61})_1\in \bigoplus_{i=0}^{m-1}\hat{A}_i$ and 
$(\ref{eq3.60})_1$, $(\ref{eq3.61})_1\in \hat{C}^{\ge 3}_{m,t,0, n-1}$, each of the right hand sides of (\ref{eq3.60}) and (\ref{eq3.61}) makes sense by the $\mathbf{k}$-linear map in (\ref{eq3.58}). It is clear that (\ref{eq3.60}) is well-defined. Using the same proof as the one which proves that (\ref{eq3.20}) is well-defined in Step 1, we get that (\ref{eq3.61}) is also well-defined. Hence, the nice  $\mathbf{k}$-linear map in (\ref{eq3.58}) can be extended to a $\mathbf{k}$-linear map in 
(\ref{eq3.59}), which is in fact nice by the same proof as the one of proving (\ref{eq3.37}).

By (\ref{eq3.55}), (\ref{eq3.57}), (\ref{eq3.58}) and (\ref{eq3.59}), (\ref{eq3.54}) holds. Using (\ref{eq3.53}) and 
(\ref{eq3.54}), in order to finish Step 2, it is enough to prove that
\begin{eqnarray}
&&\mbox{\it For $m\in\mathcal{Z}_{\ge 3}$, $m\ge t\ge 3$ and $n\in \mathcal{Z}_{\ge 1}$, a nice  $\mathbf{k}$-linear map}
\nonumber\\
\label{eq3.62}
&&\qquad\qquad\qquad\qquad\qquad\sigma : \hat{A}^{\ge 3}_{m,t,n-1}\to W \\
&&\mbox{\it can be extended to a nice  $\mathbf{k}$-linear map }\nonumber\\
\label{eq3.63}
&&\qquad\qquad\qquad\qquad\qquad\sigma : \hat{A}^{\ge 3}_{m,t,n} \to W, \\
&&\mbox{\it where $\hat{A}^{\ge 3}_{m,t,n}:=\left(\displaystyle\bigoplus_{i=0}^{m-1}\hat{A}_i\right)\bigoplus
\left(\displaystyle\bigoplus_{i=0}^{t-1}\hat{A}_m^{(i, m-i)}\right)\bigoplus
\left(\displaystyle\bigoplus_{i=0}^{n}\hat{C}^{\ge 3}_{m,t,i}\right)$.}\nonumber
\end{eqnarray}

Note that $\hat{A}^{\ge 3}_{m,t,n}=\hat{A}^{\ge 3}_{m,t,n-1}\oplus \hat{C}^{\ge 3}_{m,t,n}$ and the set $S_{m,t,n}$ defined by
$$
\left\{ \hat{x}_{j_1}\cdots\hat{x}_{j_{t-1}}\hat{x}_{j_t}\hat{q}
\hat{x}_{j_{t+1}}\cdots\hat{x}_{j_m}\left|\, \begin{array}{c}
x_{j_1}, \dots , x_{j_m}\in X,\\
ind_{\ell}(\hat{x}_{j_1}\cdots\hat{x}_{j_{t-1}}\hat{x}_{j_t}\hat{q}
\hat{x}_{j_{t+1}}\cdots\hat{x}_{j_m})=n
\end{array}\right.\right\}
$$
is a basis of $\hat{C}^{\ge 3}_{m,t,n}$. For $\hat{x}_{j_1}\cdots\hat{x}_{j_{t-1}}\hat{x}_{j_t}\hat{q}
\hat{x}_{j_{t+1}}\cdots\hat{x}_{j_m}\in S_{m,t,n}$, there exists $s$ such that 
$t-1\ge s+1>s\ge 1$ and $j_s>j_{s+1}$. Using the pair $(j_s, j_{s+1})$, we define
\begin{eqnarray*}
&&\sigma\Big(\hat{x}_{j_1}\cdots\hat{x}_{j_{s-1}}\cdot\hat{x}_{j_s} \hat{x}_{j_{s+1}}\cdot\hat{x}_{j_{s+2}}
\cdots \hat{x}_{j_{t-1}}\hat{x}_{j_t}\hat{q}\hat{x}_{j_{t+1}}\cdots\hat{x}_{j_m}\Big)\\
&:=&\sigma\Big(\underbrace{\hat{x}_{j_1}\cdots\hat{x}_{j_{s-1}}\cdot [\hat{x}_{j_s}, \hat{x}_{j_{s+1}}]\,\hat{}\;\cdot\hat{x}_{j_{s+2}}
\cdots \hat{x}_{j_{t-1}}\hat{x}_{j_t}\hat{q}\hat{x}_{j_{t+1}}\cdots\hat{x}_{j_m}}_{(\ref{eq3.64})_1}\Big)+\\
&&+\sigma\Big(\underbrace{\hat{x}_{j_1}\cdots\hat{x}_{j_{s-1}}\cdot(\hat{x}_{j_{s+1}}\hat{x}_{j_s})\cdot\hat{x}_{j_{s+2}}
\cdots \hat{x}_{j_{t-1}}\hat{x}_{j_t}\hat{q}\hat{x}_{j_{t+1}}\cdots\hat{x}_{j_m}}_{(\ref{eq3.64})_2}\Big)+\\
&&+\sigma\Big(\underbrace{\hat{x}_{j_1}\cdots\hat{x}_{j_{s-1}}\cdot(\hat{x}_{j_s}\hat{x}_{j_{s+1}}\hat{q})\cdot\hat{x}_{j_{s+2}}
\cdots \hat{x}_{j_{t-1}}\hat{x}_{j_t}\hat{x}_{j_{t+1}}\cdots\hat{x}_{j_m}}_{(\ref{eq3.64})_3}\Big)+\\
&&-\sigma\Big(\underbrace{\hat{x}_{j_1}\cdots\hat{x}_{j_{s-1}}\cdot(\hat{x}_{j_{s+1}}\hat{x}_{j_s}\hat{q})\cdot\hat{x}_{j_{s+2}}
\cdots \hat{x}_{j_{t-1}}\hat{x}_{j_t}\hat{x}_{j_{t+1}}\cdots\hat{x}_{j_m}}_{(\ref{eq3.64})_4}\Big)+\\
&&-k\sigma\Big(\underbrace{\hat{x}_{j_1}\cdots\hat{x}_{j_{s-1}}\cdot(\hat{x}_{j_s}\hat{q}\hat{x}_{j_{s+1}})\cdot\hat{x}_{j_{s+2}}
\cdots \hat{x}_{j_{t-1}}\hat{x}_{j_t}\hat{x}_{j_{t+1}}\cdots\hat{x}_{j_m}}_{(\ref{eq3.64})_5}\Big)+
\end{eqnarray*} 
\begin{equation}\label{eq3.64}
+k\sigma\Big(\underbrace{\hat{x}_{j_1}\cdots\hat{x}_{j_{s-1}}\cdot(\hat{x}_{j_{s+1}}\hat{q}\hat{x}_{j_s})\cdot\hat{x}_{j_{s+2}}
\cdots \hat{x}_{j_{t-1}}\hat{x}_{j_t}\hat{x}_{j_{t+1}}\cdots\hat{x}_{j_m}}_{(\ref{eq3.64})_6}\Big).
\end{equation} 

Since $deg_X\left((\ref{eq3.64})_1\right)=m-1$,  $ind_{\ell}\left((\ref{eq3.64})_2\right)=n-1$,  
$deg_{\hat{q}}\left((\ref{eq3.64})_3\right)=deg_{\hat{q}}\left((\ref{eq3.64})_4\right)=s+1<t$ and
$deg_{\hat{q}}\left((\ref{eq3.64})_5\right)=deg_{\hat{q}}\left((\ref{eq3.64})_6\right)=s<t$, each of the six terms on the right hand side of 
(\ref{eq3.64}) makes sense by (\ref{eq3.62}). 

\medskip
One can prove that the right hand side of 
(\ref{eq3.64}) is well-defined.
Hence, we have extended the nice $\mathbf{k}$-linear map in (\ref{eq3.62}) to a $\mathbf{k}$-linear map in (\ref{eq3.63}) satisfying (\ref{eq3.64}). The remaining part of {\it Step 2} is to prove that the extended $\mathbf{k}$-linear map in (\ref{eq3.63}) satisfying (\ref{eq3.64}) is nice.

\medskip
By (\ref{eq3.2}), (\ref{eq3.3}) and (\ref{eq3.4}), if a  $(u, v)^{[6-th]}$-generator $H_{u,v}$ satisfies
$$
H_{u,v}\in \hat{A}^{\ge 3}_{m,t,n}=\left(\displaystyle\bigoplus_{i=0}^{m-1}\hat{A}_i\right)\bigoplus
\left(\displaystyle\bigoplus_{i=0}^{t-1}\hat{A}_m^{(i, m-i)}\right)\bigoplus
\left(\displaystyle\bigoplus_{i=0}^{n}\hat{C}^{\ge 3}_{m,t,i}\right)
$$
and
$$
H_{u,v}\not\in \hat{A}^{\ge 3}_{m,t,n-1}=\left(\displaystyle\bigoplus_{i=0}^{m-1}\hat{A}_i\right)\bigoplus
\left(\displaystyle\bigoplus_{i=0}^{t-1}\hat{A}_m^{(i, m-i)}\right)\bigoplus
\left(\displaystyle\bigoplus_{i=0}^{n-1}\hat{C}^{\ge 3}_{m,t,i}\right),
$$
then $(u, v)=(1,0)$ or $(0, 1)$. Clearly, we have
\begin{equation}\label{eq3.81}
H_{0,1}\in \hat{A}^{\ge 3}_{m,t,n}\setminus\hat{A}^{\ge 3}_{m,t,n-1}\Longrightarrow \sigma(H_{0,1})=0
\end{equation}
by (\ref{eq3.64}).

If $H_{1, 0}\in \hat{A}^{\ge 3}_{m,t,n}\setminus\hat{A}^{\ge 3}_{m,t,n-1}$, then $H_{1, 0}$ can be written as 
$$
\hat{x}_{j_1}\cdots\hat{x}_{j_{r-1}}(\hat{x}_{j_r}\hat{x}_{j_{r+1}})\hat{x}_{j_{r+2}}\cdots
\hat{x}_{j_t}\hat{q}\hat{x}_{j_{t+1}}\cdots\hat{x}_{j_{s-1}}
[x_{j_s}, x_{j_{s+1}}]\,\hat{}\;\hat{x}_{j_{s+2}}\cdots\hat{x}_{j_m}+
$$
$$
-k\hat{x}_{j_1}\cdots\hat{x}_{j_{r-1}}(\hat{x}_{j_r}\hat{x}_{j_{r+1}})\hat{x}_{j_{r+2}}\cdots
\hat{x}_{j_t}\hat{q}\hat{x}_{j_{t+1}}\cdots\hat{x}_{j_{s-1}}
(\hat{x}_{j_s}\hat{x}_{j_{s+1}})\hat{x}_{j_{s+2}}\cdots\hat{x}_{j_m}+
$$
$$
+k\hat{x}_{j_1}\cdots\hat{x}_{j_{r-1}}(\hat{x}_{j_r}\hat{x}_{j_{r+1}})\hat{x}_{j_{r+2}}\cdots
\hat{x}_{j_t}\hat{q}\hat{x}_{j_{t+1}}\cdots\hat{x}_{j_{s-1}}
(\hat{x}_{j_{s+1}}\hat{x}_{j_s})\hat{x}_{j_{s+2}}\cdots\hat{x}_{j_m}
$$
by (\ref{eq3.3}), where $0< r<r+1< t$, $j_r>j_{r+1}$ and each of the three terms above has the left index $n$. Applying (\ref{eq3.64}) to each of the three terms above, we get
$$
\sigma(H_{1, 0})=
$$
$$\footnotesize
\underbrace{\sigma\Big(\hat{x}_{j_1}\cdots\hat{x}_{j_{r-1}}\cdot [x_{j_r}, x_{j_{r+1}}]\,\hat{}\;\cdot\hat{x}_{j_{r+2}}\cdots
\hat{x}_{j_t}\hat{q}\hat{x}_{j_{t+1}}\cdots\hat{x}_{j_{s-1}}
[x_{j_s}, x_{j_{s+1}}]\,\hat{}\;\hat{x}_{j_{s+2}}\cdots\hat{x}_{j_m}\Big)}_{(\ref{eq3.82})_{1,1}}+
$$
$$\footnotesize
+\underbrace{\sigma\Big(\hat{x}_{j_1}\cdots\hat{x}_{j_{r-1}}(\hat{x}_{j_{r+1}}\hat{x}_{j_r})\hat{x}_{j_{r+2}}\cdots
\hat{x}_{j_t}\hat{q}\hat{x}_{j_{t+1}}\cdots\hat{x}_{j_{s-1}}
[x_{j_s}, x_{j_{s+1}}]\,\hat{}\;\hat{x}_{j_{s+2}}\cdots\hat{x}_{j_m}\Big)}_{(\ref{eq3.82})_{1,2}}+
$$
$$\footnotesize
+\underbrace{\sigma\Big(\hat{x}_{j_1}\cdots\hat{x}_{j_{r-1}}(\hat{x}_{j_r}\hat{x}_{j_{r+1}}\hat{q})\hat{x}_{j_{r+2}}\cdots
\hat{x}_{j_t}\hat{q}\hat{x}_{j_{t+1}}\cdots\hat{x}_{j_{s-1}}
[x_{j_s}, x_{j_{s+1}}]\,\hat{}\;\hat{x}_{j_{s+2}}\cdots\hat{x}_{j_m}\Big)}_{(\ref{eq3.82})_{1,3}}+
$$
$$\footnotesize
-\underbrace{\sigma\Big(\hat{x}_{j_1}\cdots\hat{x}_{j_{r-1}}(\hat{x}_{j_{r+1}}\hat{x}_{j_r}\hat{q})\hat{x}_{j_{r+2}}\cdots
\hat{x}_{j_t}\hat{x}_{j_{t+1}}\cdots\hat{x}_{j_{s-1}}
[x_{j_s}, x_{j_{s+1}}]\,\hat{}\;\hat{x}_{j_{s+2}}\cdots\hat{x}_{j_m}\Big)}_{(\ref{eq3.82})_{1,4}}+
$$
$$\footnotesize
-k\underbrace{\sigma\Big(\hat{x}_{j_1}\cdots\hat{x}_{j_{r-1}}(\hat{x}_{j_r}\hat{q}\hat{x}_{j_{r+1}})\hat{x}_{j_{r+2}}\cdots
\hat{x}_{j_t}\hat{x}_{j_{t+1}}\cdots\hat{x}_{j_{s-1}}
[x_{j_s}, x_{j_{s+1}}]\,\hat{}\;\hat{x}_{j_{s+2}}\cdots\hat{x}_{j_m}\Big)}_{(\ref{eq3.82})_{1,5}}+
$$
$$\footnotesize
+k\underbrace{\sigma\Big(\hat{x}_{j_1}\cdots\hat{x}_{j_{r-1}}(\hat{x}_{j_{r+1}}\hat{q}\hat{x}_{j_r})\hat{x}_{j_{r+2}}\cdots
\hat{x}_{j_t}\hat{x}_{j_{t+1}}\cdots\hat{x}_{j_{s-1}}
[x_{j_s}, x_{j_{s+1}}]\,\hat{}\;\hat{x}_{j_{s+2}}\cdots\hat{x}_{j_m}\Big)}_{(\ref{eq3.82})_{1,6}}+
$$
$$\footnotesize
-k\underbrace{\sigma\Big(\hat{x}_{j_1}\cdots\hat{x}_{j_{r-1}}\cdot [x_{j_r}, x_{j_{r+1}}]\,\hat{}\;\cdot \hat{x}_{j_{r+2}}\cdots
\hat{x}_{j_t}\hat{q}\hat{x}_{j_{t+1}}\cdots\hat{x}_{j_{s-1}}
(\hat{x}_{j_s}\hat{x}_{j_{s+1}})\hat{x}_{j_{s+2}}\cdots\hat{x}_{j_m}\Big)}_{(\ref{eq3.82})_{2,1}}+
$$
$$\footnotesize
-k\underbrace{\sigma\Big(\hat{x}_{j_1}\cdots\hat{x}_{j_{r-1}}(\hat{x}_{j_{r+1}}\hat{x}_{j_r})\hat{x}_{j_{r+2}}\cdots
\hat{x}_{j_t}\hat{q}\hat{x}_{j_{t+1}}\cdots\hat{x}_{j_{s-1}}
(\hat{x}_{j_s}\hat{x}_{j_{s+1}})\hat{x}_{j_{s+2}}\cdots\hat{x}_{j_m}\Big)}_{(\ref{eq3.82})_{2,2}}+
$$
$$\footnotesize
-k\underbrace{\sigma\Big(\hat{x}_{j_1}\cdots\hat{x}_{j_{r-1}}(\hat{x}_{j_r}\hat{x}_{j_{r+1}}\hat{q})\hat{x}_{j_{r+2}}\cdots
\hat{x}_{j_t}\hat{x}_{j_{t+1}}\cdots\hat{x}_{j_{s-1}}
(\hat{x}_{j_s}\hat{x}_{j_{s+1}})\hat{x}_{j_{s+2}}\cdots\hat{x}_{j_m}\Big)}_{(\ref{eq3.82})_{2,3}}+
$$
$$\footnotesize
+k\underbrace{\sigma\Big(\hat{x}_{j_1}\cdots\hat{x}_{j_{r-1}}(\hat{x}_{j_{r+1}}\hat{x}_{j_r}\hat{q})\hat{x}_{j_{r+2}}\cdots
\hat{x}_{j_t}\hat{x}_{j_{t+1}}\cdots\hat{x}_{j_{s-1}}
(\hat{x}_{j_s}\hat{x}_{j_{s+1}})\hat{x}_{j_{s+2}}\cdots\hat{x}_{j_m}\Big)}_{(\ref{eq3.82})_{2,4}}+
$$
$$\footnotesize
+k^2\underbrace{\sigma\Big(\hat{x}_{j_1}\cdots\hat{x}_{j_{r-1}}(\hat{x}_{j_r}\hat{q}\hat{x}_{j_{r+1}})\hat{x}_{j_{r+2}}\cdots
\hat{x}_{j_t}\hat{x}_{j_{t+1}}\cdots\hat{x}_{j_{s-1}}
(\hat{x}_{j_s}\hat{x}_{j_{s+1}})\hat{x}_{j_{s+2}}\cdots\hat{x}_{j_m}\Big)}_{(\ref{eq3.82})_{2,5}}+
$$
$$\footnotesize
-k^2\underbrace{\sigma\Big(\hat{x}_{j_1}\cdots\hat{x}_{j_{r-1}}(\hat{x}_{j_{r+1}}\hat{q}\hat{x}_{j_r})\hat{x}_{j_{r+2}}\cdots
\hat{x}_{j_t}\hat{x}_{j_{t+1}}\cdots\hat{x}_{j_{s-1}}
(\hat{x}_{j_s}\hat{x}_{j_{s+1}})\hat{x}_{j_{s+2}}\cdots\hat{x}_{j_m}\Big)}_{(\ref{eq3.82})_{2,6}}+
$$
$$\footnotesize
+k\underbrace{\sigma\Big(\hat{x}_{j_1}\cdots\hat{x}_{j_{r-1}}\cdot [x_{j_r}, x_{j_{r+1}}]\,\hat{}\;\cdot \hat{x}_{j_{r+2}}\cdots
\hat{x}_{j_t}\hat{q}\hat{x}_{j_{t+1}}\cdots\hat{x}_{j_{s-1}}
(\hat{x}_{j_{s+1}}\hat{x}_{j_s})\hat{x}_{j_{s+2}}\cdots\hat{x}_{j_m}\Big)}_{(\ref{eq3.82})_{3,1}}+
$$
$$\footnotesize
+k\underbrace{\sigma\Big(\hat{x}_{j_1}\cdots\hat{x}_{j_{r-1}}(\hat{x}_{j_{r+1}}\hat{x}_{j_r})\hat{x}_{j_{r+2}}\cdots
\hat{x}_{j_t}\hat{q}\hat{x}_{j_{t+1}}\cdots\hat{x}_{j_{s-1}}
(\hat{x}_{j_{s+1}}\hat{x}_{j_s})\hat{x}_{j_{s+2}}\cdots\hat{x}_{j_m}\Big)}_{(\ref{eq3.82})_{3,2}}+
$$
$$\footnotesize
+k\underbrace{\sigma\Big(\hat{x}_{j_1}\cdots\hat{x}_{j_{r-1}}(\hat{x}_{j_r}\hat{x}_{j_{r+1}}\hat{q})\hat{x}_{j_{r+2}}\cdots
\hat{x}_{j_t}\hat{x}_{j_{t+1}}\cdots\hat{x}_{j_{s-1}}
(\hat{x}_{j_{s+1}}\hat{x}_{j_s})\hat{x}_{j_{s+2}}\cdots\hat{x}_{j_m}\Big)}_{(\ref{eq3.82})_{3,3}}+
$$
$$\footnotesize
-k\underbrace{\sigma\Big(\hat{x}_{j_1}\cdots\hat{x}_{j_{r-1}}(\hat{x}_{j_{r+1}}\hat{x}_{j_r}\hat{q})\hat{x}_{j_{r+2}}\cdots
\hat{x}_{j_t}\hat{x}_{j_{t+1}}\cdots\hat{x}_{j_{s-1}}
(\hat{x}_{j_{s+1}}\hat{x}_{j_s})\hat{x}_{j_{s+2}}\cdots\hat{x}_{j_m}\Big)}_{(\ref{eq3.82})_{3,4}}+
$$
$$\footnotesize
-k^2\underbrace{\sigma\Big(\hat{x}_{j_1}\cdots\hat{x}_{j_{r-1}}(\hat{x}_{j_r}\hat{q}\hat{x}_{j_{r+1}})\hat{x}_{j_{r+2}}\cdots
\hat{x}_{j_t}\hat{x}_{j_{t+1}}\cdots\hat{x}_{j_{s-1}}
(\hat{x}_{j_{s+1}}\hat{x}_{j_s})\hat{x}_{j_{s+2}}\cdots\hat{x}_{j_m}\Big)}_{(\ref{eq3.82})_{3,5}}+
$$
\begin{equation}\label{eq3.82}
\footnotesize
+k^2\underbrace{\sigma\Big(\hat{x}_{j_1}\cdots\hat{x}_{j_{r-1}}(\hat{x}_{j_{r+1}}\hat{q}\hat{x}_{j_r})\hat{x}_{j_{r+2}}\cdots
\hat{x}_{j_t}\cdots\hat{x}_{j_{s-1}}
(\hat{x}_{j_{s+1}}\hat{x}_{j_s})\hat{x}_{j_{s+2}}\cdots\hat{x}_{j_m}\Big)}_{(\ref{eq3.82})_{3,6}}.
\end{equation}

The sum
$(\ref{eq3.82})_{1,i}-(\ref{eq3.82})_{2,i}+(\ref{eq3.82})_{3,i}$ for $6\ge i\ge 1$
is the image of a  $(1, 0)^{[6-th]}$-generator in $\hat{A}^{\ge 3}_{m,t,n-1}$ under the nice $\mathbf{k}$-linear map in (\ref{eq3.62}). Hence, we have
$$
(\ref{eq3.82})_{1,i}-(\ref{eq3.82})_{2,i}+(\ref{eq3.82})_{3,i}=0\quad\mbox{ for $6\ge i\ge 1$,}
$$
which implies that
\begin{equation}\label{eq3.83}
\sigma(H_{1, 0})=0 \quad\mbox{$H_{1, 0}\in \hat{A}^{\ge 3}_{m,t,n}\setminus\hat{A}^{\ge 3}_{m,t,n-1}$}
\end{equation}
by (\ref{eq3.82}). It follows from (\ref{eq3.81}) and (\ref{eq3.83}) that the extended $\mathbf{k}$-linear map in (\ref{eq3.63}) satisfying (\ref{eq3.64}) is nice.

\medskip
This completes {\it Step 2}.

\medskip
\subsection{Step 3}

The goal of Step 3 is to prove the following fact.
\begin{center}
{\it  {\bf Fact 4:} For $m\in\mathcal{Z}_{\ge 2}$, a nice $\mathbf{k}$-linear map}\\
\begin{equation}\label{eq3.84}
\sigma : \hat{A}_{m-1, m, m}:=\left(\displaystyle\bigoplus_{i=0}^{m-1}\hat{A}_i\right)\bigoplus
\left(\displaystyle\bigoplus_{i=0}^{m}\hat{A}_m^{(i, m-i)}\right)\to W
\end{equation}\\
{\it can be extended to a nice $\mathbf{k}$-linear map}\\
\begin{equation}\label{eq3.85}
\sigma : \hat{A}_{m-1, m, m+1}:=\left(\displaystyle\bigoplus_{i=0}^{m-1}\hat{A}_i\right)\bigoplus
\left(\displaystyle\bigoplus_{i=0}^{m+1}\hat{A}_m^{(i, m-i)}
\right)=\displaystyle\bigoplus_{i=0}^{m}\hat{A}_i\to W.
\end{equation}
\end{center}

Recall that
\begin{equation}\label{eq3.86}
\hat{A}_m^{(m+1, -1)}=\bigoplus_{i\in\mathcal{Z}_{\ge 0}}\hat{A}_{m,i}^{(m+1, -1)},
\end{equation}
where $\hat{A}_{m,i}^{(m+1, -1)}$ is the subspace spanned by all monomials whose $X$-degrees are $m$, whose $\hat{q}$-degrees are 
$(m+1, -1)$, and whose right indices are $i$.

\medskip
First, we define
\begin{equation}\label{eq3.87}
\sigma\left(\hat{x}_{j_1}\cdots\hat{x}_{j_m}\right):=\mathring{x}_{j_1}\cdots\mathring{x}_{j_m}\quad\mbox{for $j_m\ge \cdots \ge j_1$.}
\end{equation}
Since 
$$\left\{\, \hat{x}_{j_1}\cdots\hat{x}_{j_m}\,|\, x_{j_1}, \cdots , x_{j_m}\in X\quad\mbox{and}\quad j_m\ge \cdots \ge j_1\,\right\}$$
is a basis of $\hat{A}_{m,0}^{(m+1, -1)}$, (\ref{eq3.87}) extends the nice $\mathbf{k}$-linear map in (\ref{eq3.84}) to a 
$\mathbf{k}$-linear map
\begin{equation}\label{eq3.88}
\sigma: \hat{A}_{m-1, m, m}\bigoplus \hat{A}_{m,0}^{(m+1, -1)}\to W.
\end{equation}

All $(u, v)^{[6-th]}$-generator in $\hat{A}_{m-1, m, m}\bigoplus \hat{A}_{m,0}^{(m+1, -1)}$ must be in $\hat{A}_{m-1, m, m}$ by 
(\ref{eq3.2}), (\ref{eq3.3}) and (\ref{eq3.4}). Hence, the $\mathbf{k}$-linear map in (\ref{eq3.88}) is nice.

\medskip
We now prove that
\begin{center}
{\it  {\bf Fact 5:} For $n\in\mathcal{Z}_{\ge 1}$, a nice $\mathbf{k}$-linear map}\\
\begin{equation}\label{eq3.89}
\sigma : \hat{A}_{m-1, m, m}\bigoplus
\left(\displaystyle\bigoplus_{i=0}^{n-1}\hat{A}_{m,i}^{(m+1, -1)}\right)\to W
\end{equation}\\
{\it can be extended to a nice $\mathbf{k}$-linear map}\\
\begin{equation}\label{eq3.90}
\sigma : \hat{A}_{m-1, m, m}\bigoplus
\left(\displaystyle\bigoplus_{i=0}^{n}\hat{A}_{m,i}^{(m+1, -1)}\right)\to W.
\end{equation}
\end{center}

The set
$$
S_n:=\{\, \hat{x}_{j_1}\cdots\hat{x}_{j_m} \,|\, x_{j_1}, \cdots , x_{j_m}\in X\quad\mbox{and}\quad 
ind_r(\hat{x}_{j_1}\cdots\hat{x}_{j_m})=n\,\}
$$
is a basis of $\hat{A}_{m,n}^{(m+1, -1)}$. For any $\hat{x}_{j_1}\cdots\hat{x}_{j_m}\in S_n$ with $n\in\mathcal{Z}_{\ge 1}$, there 
exists $s$ such that $m\ge s+1>s\ge 1$ and $j_s>j_{s+1}$. Using the pair $(j_s, j_{s+1})$, we define
\begin{eqnarray}\label{eq3.91}
&&\sigma\Big(\hat{x}_{j_1}\cdots\hat{x}_{j_{s-1}}\cdot\hat{x}_{j_s}\hat{x}_{j_{s+1}}\cdot\hat{x}_{j_{s+2}}
\cdots\hat{x}_{j_m} \Big)\nonumber\\
&=&\sigma\Big(\hat{x}_{j_1}\cdots\hat{x}_{j_{s-1}}\cdot [x_{s}, x_{j_{s+1}}]\,\hat{}\cdot\hat{x}_{j_{s+2}}
\cdots\hat{x}_{j_m} \Big)+\nonumber\\
&&+\sigma\Big(\hat{x}_{j_1}\cdots\hat{x}_{j_{s-1}}(\hat{x}_{j_{s+1}}\hat{x}_{j_s})\hat{x}_{j_{s+2}}
\cdots\hat{x}_{j_m} \Big)+\nonumber\\
&&+\sigma\Big(\hat{x}_{j_1}\cdots\hat{x}_{j_{s-1}}(\hat{x}_{j_s}\hat{x}_{j_{s+1}}\hat{q})\hat{x}_{j_{s+2}}
\cdots\hat{x}_{j_m} \Big)+\nonumber\\
&&-\sigma\Big(\hat{x}_{j_1}\cdots\hat{x}_{j_{s-1}}(\hat{x}_{j_{s+1}}\hat{x}_{j_s}\hat{q})\hat{x}_{j_{s+2}}
\cdots\hat{x}_{j_m} \Big)+\nonumber\\
&&-k\sigma\Big(\hat{x}_{j_1}\cdots\hat{x}_{j_{s-1}}(\hat{x}_{j_s}\hat{q}\hat{x}_{j_{s+1}})\hat{x}_{j_{s+2}}
\cdots\hat{x}_{j_m} \Big)+\nonumber\\
&&+k\sigma\Big(\hat{x}_{j_1}\cdots\hat{x}_{j_{s-1}}(\hat{x}_{j_{s+1}}\hat{q}\hat{x}_{j_s})\hat{x}_{j_{s+2}}
\cdots\hat{x}_{j_m} \Big).
\end{eqnarray}
Using the $\mathbf{k}$-linear map in (\ref{eq3.89}), each term on the right hand side of (\ref{eq3.91}) makes sense. By the same proof as the one of proving that (\ref{eq3.64}) is well-defined in {\it Step 2}, we know that (\ref{eq3.91}) is well-defined. Thus, we get a 
$\mathbf{k}$-linear map in (\ref{eq3.90}).

If a  $(u, v)^{[6-th]}$-generator $H_{u,v}$ satisfies
\begin{equation}\label{eq3.92}
H_{u,v}\in \hat{A}_{m-1, m, m}\bigoplus
\left(\displaystyle\bigoplus_{i=0}^{n}\hat{A}_{m,i}^{(m+1, -1)}\right)
\end{equation}
and
\begin{equation}\label{eq3.93}
H_{u,v}\not\in \hat{A}_{m-1, m, m}\bigoplus
\left(\displaystyle\bigoplus_{i=0}^{n-1}\hat{A}_{m,i}^{(m+1, -1)}\right),
\end{equation}
then $(u, v)=(0,0)$ by (\ref{eq3.2}), (\ref{eq3.3}) and (\ref{eq3.4}). It is also clear that if $H_{00}$ satisfies (\ref{eq3.92}) and 
(\ref{eq3.93}), then $\sigma(H_{00})=0$ by (\ref{eq3.91}). Hence, the $\mathbf{k}$-linear map in (\ref{eq3.90}), which extends the $\mathbf{k}$-linear map in (\ref{eq3.89}) and satisfies (\ref{eq3.91}), is nice. This fact and (\ref{eq3.86}) imply that the $\mathbf{k}$-linear map in (\ref{eq3.84}) can be extended to a nice $\mathbf{k}$-linear map in (\ref{eq3.85}). Hence, {\it Step 3} is done.

\medskip
Using the facts proved in the three steps above, we can now prove  Proposition~\ref{pr3.2}. 

\medskip
By (\ref{eq3.8}), there exists a nice $\mathbf{k}$-linear map 
\begin{equation}\label{eq3.94}
\sigma : \hat{A}_0\oplus \hat{A}_1\to W,
\end{equation}
which can be extended to a nice    
$\mathbf{k}$-linear map
\begin{equation}\label{eq3.95}
\sigma: \left(\displaystyle\bigoplus_{i=0}^{m-1}\hat{A}_i\right)\bigoplus
\left(\displaystyle\bigoplus_{i=0}^{m}\hat{A}_m^{(i, m-i)}
\right)\to W \quad\mbox{with $m\in\mathcal{Z}_{\ge 2}$}
\end{equation}
by {\it Step 1} and {\it Step 2}. It follows from  {\it Step 3} that the nice $\mathbf{k}$-linear map in (\ref{eq3.95}) can be extended to a nice  $\mathbf{k}$-linear map
$$
\sigma: \displaystyle\bigoplus_{i=0}^{m}\hat{A}_i=\left(\displaystyle\bigoplus_{i=0}^{m-1}\hat{A}_i\right)\bigoplus
\left(\displaystyle\bigoplus_{i=0}^{m+1}\hat{A}_m^{(i, m-i)}
\right)\to W \quad\mbox{with $m\in\mathcal{Z}_{\ge 2}$.}
$$
In other words, the nice $\mathbf{k}$-linear map in (\ref{eq3.94})  can be extended to a nice  $\mathbf{k}$-linear map
$\sigma: (\hat{A}, \hat{q})=\displaystyle\bigoplus_{m=0}^{\infty}\hat{A}_m\to W$. This proves that Proposition~\ref{pr3.2} is true.

\bigskip
The extended$^{6-th}$ P-B-W Theorem is the following

\begin{proposition}\label{pr3.3} {\bf (The Extended$^{6-th}$ P-B-W Theorem)}  Let $\mathcal{L}$ be a  Lie algebra with a basis  
$X:=\{\, x_j \,|\, j\in J\,\}$. If
$\Big(\mathcal{U}:=\displaystyle\frac{\hat{A}}{R}, \quad \bar{q}:=\hat{q}+R\Big)$ is the enveloping$^{6-th}$ algebra for the Lie algebra 
$\mathcal{L}$, then
the following set of cosets of  model monomials
$$
\hat{T}:=\left\{\begin{array}{c}
\hat{q}\hat{x}_{j_1}\cdots\hat{x}_{j_m}+R,\\
\hat{x}_{i_1}\cdots \hat{x}_{i_t}\hat{x}_{j_0}\hat{q}\hat{x}_{j_1}\cdots\hat{x}_{j_m}+R, \\
\hat{x}_{j_1}\cdots\hat{x}_{j_m}+R \\\end{array}
 \,\left|\, 
\begin{array}{c}x_{i_1}, \cdots , x_{i_t}, x_{j_0}, x_{j_1}, \cdots , x_{j_m}\in X,\\
i_t\geq\cdots \geq i_1,\\ 
j_m\ge \cdots\ge j_1\geq j_0,\\
t, m\in\mathcal{Z}_{\ge 0}
\end{array}\right.\right\}.
$$
is a basis for the enveloping$^{6-th}$ algebra $(\mathcal{U}, \bar{q}))$.
\end{proposition}

\medskip
\noindent
{\bf Proof} By Proposition~\ref{pr3.2}, there is a nice $\mathbf{k}$-linear map $\sigma : (\hat{A}, \hat{q})\to W$. Thus, $\sigma(R)=0$. Hence, the nice $\mathbf{k}$-linear map $\sigma : (\hat{A}, \hat{q})\to W$ induces a $\mathbf{k}$-linear map: 
$\bar{\sigma}: \mathcal{U}=\displaystyle\frac{\hat{A}}{R}\to W$. On one hand, the image of the cosets of  model monomials in $\hat{T}$ are linearly independent by (\ref{eq3.5}). On the other hand, the cosets of  model monomials in $\hat{T}$ also span the vector space 
$\mathcal{U}=\displaystyle\frac{\hat{A}}{R}$. Hence, Proposition~\ref{pr3.3} holds.

\hfill\raisebox{1mm}{\framebox[2mm]{}}

\section{Hopf-like Algebras}

In this section, we discuss the Hopf-like algebra structure on the  enveloping$^{6-th}$ algebra
of a Lie algebra.

\medskip
Let $(A, q)$ be an invariant algebra induced by the idempotent $q$. We define a product $\circ$ in $\mathring{A}:=A$ by
\begin{equation}\label{eq5.1}
x\circ y:=yx\qquad \mbox{for $x, y\in \mathring{A}$.}
\end{equation}

Since
$$
(1-q)\circ (1-q)=(1-q)(1-q)=1-2q+q^2=1-q
$$
and
$$
(1-q)\circ x\circ (1-q)=(1-q)x(1-q)=x(1-q)=(1-q)\circ x\qquad \mbox{for $x\in \mathring{A}$,}
$$
$(\mathring{A}, 1-q)$ is an invariant algebra induced by the idempotent $1-q$, which is called the {\bf opposite ( right) invariant algebra} of the invariant algebra $(A, q)$. An {\bf invariant anti-homomorphism} $f$ of an invariant algebra $(A, q)$ is an invariant homomorphism
$f: (A, q)\to (\mathring{A}, 1-q)$.
By (\ref{eq1.4}), the square bracket $[\, ,\,]^{\circ}_{6, k}$ which defines the Lie algebra 
$Lie\Big ((\mathring{A}, 1-q), [\, ,\,]^{\circ}_{6, k}\Big)$ is given by 
$$
[x , y]^{\circ}_{6, k}:=x\circ y-y\circ x-x\circ y\circ (1-q)+y\circ x\circ (1-q)
+kx\circ (1-q)\circ y-ky\circ (1-q)\circ x
$$
or
\begin{equation}\label{eq5.2}
[x , y]^{\circ}_{6, k}=qyx-qxy+kyx-kyqx-kxy+kxqy,
\end{equation}
where $x, y\in \mathring{A}=A$.

\begin{proposition}\label{pr5.1} Let $\Big(U^{6-th}(L), q\Big)$ be the 
enveloping$^{6-th}$ algebra of a Lie algebra $L$. There exists a unique invariant anti-homomorphism $S$ of 
\newline $\Big(U^{6-th}(L), q\Big)$ such that $S(q)=1-q$ and
\begin{equation}\label{eq5.3}
S(x)=-\frac1k x-kqx+\frac1k xq\qquad\mbox{for $x\in L$.}
\end{equation}
\end{proposition}

\medskip
\noindent
{\bf Proof} It follows from the universal property of the enveloping$^{6-th}$ algebra.

\hfill\raisebox{1mm}{\framebox[2mm]{}}

\medskip
Let $C$ be a vector space over a field $\mathbf{k}$. The canonical map $c\mapsto 1\otimes c$ is denoted by  $C\rightarrow\mathbf{k}\otimes C$, and the canonical map 
$c\mapsto c\otimes 1$ is denoted by $C\rightarrow C\otimes \mathbf{k}$  respectively, where $c\in C$.

\medskip
We now introduce coalgebras with $\sigma$-counits in the following

\begin{definition}\label{def5.1} A vector space $C$ over a field $\mathbf{k}$ is called a {\bf coalgebra with $\sigma$-counit} if there exist three linear maps $\Delta : C\to C\otimes C$,
$\varepsilon : C\to \mathbf{k}$ and $\sigma : C\to C$ such that the diagram 
$$
\begin{array}{ccc}
C&\stackrel{\Delta}{\longrightarrow}&C\otimes C\\
\Delta\Bigg\downarrow&&\Bigg\downarrow id\otimes\Delta\\
C\otimes C&\stackrel{\Delta\otimes id}{\longrightarrow}&C\otimes C\otimes C
\end{array}
$$
and the diagram
$$
\begin{array}{ccccc}
\mathbf{k}\otimes C&\longleftarrow&C&\longrightarrow
&C\otimes \mathbf{k}\\
\varepsilon\otimes id\Bigg\uparrow&&\Bigg\uparrow\sigma&&\Bigg\uparrow id\otimes \varepsilon\\
C\otimes C&\stackrel{\Delta}{\longleftarrow}&C&\stackrel{\Delta}{\longrightarrow}
&C\otimes C
\end{array}
$$
are commutative. The maps $\Delta$ and $\varepsilon$ are called the 
{\bf comultiplication} and the {\bf $\sigma$-counit} respectively.
\end{definition}

\medskip
Clearly, an ordinary counit is an $id$-counit. Hence, an ordinary coalgebra is a coalgebra with $id$-counit.

\medskip
We now construct an  algebra endomorphism of an invariant algebra.

\begin{proposition}\label{pr5.2} Let $(A, q)$ be an invariant algebra. If 
$\sigma : (A, q)\to (A, q)$ is a linear map defined by
\begin{equation}\label{eq5.7}
\sigma (a):=a+qa-aq\qquad\mbox{for $a\in (A, q)$,}
\end{equation}
then $\sigma$ has the following properties:
\begin{enumerate}
\item $\sigma ^2=\sigma$;
\item $\sigma$ is an invariant homomorphism;
\item $\sigma -id$ is a derivation of $(A, q)$.
\end{enumerate}
\end{proposition}

\medskip
\noindent
{\bf Proof} Let $a$ and $b$ be two elements of $(A, q)$.
Since 
\begin{eqnarray*}
&&\sigma ^2(a)=\sigma (a+qa-aq)\\
&=&a+qa-aq+q(a+qa-aq)-(a+qa-aq)q\\
&=&a+qa-aq+qa-qa=a+qa-aq=\sigma (a),
\end{eqnarray*}
which proves the property 1. 

Similarly, both the property 2 and the property 3 hold.

\hfill\raisebox{1mm}{\framebox[2mm]{}}

\bigskip
The following proposition proves that the enveloping$^{6-th}$ algebra of a Lie algebra is a 
coalgebra with $\sigma$-counit.

\begin{proposition}\label{pr5.3} If $\Big(U^{6-th}(L), q\Big)$ is the 
enveloping$^{6-th}$ algebra of a Lie algebra $L$, then there exist two algebra homomorphisms 
$\varepsilon :  U^{6-th}(L)\to\mathbf{k}$ and 
$\Delta : U^{6-th}(L)\to U^{6-th}(L)\otimes U^{6-th}(L)$ such that 
\begin{equation}\label{eq5.8}
\Delta (q)=q\otimes q ,
\end{equation}
\begin{eqnarray}\label{eq5.9}
\Delta (x)&=&(x+kqx-xq)\otimes 1+1\otimes (x+kqx-xq)+\nonumber\\
&& +(1-k)qx\otimes q+(1-k)q\otimes qx \qquad\mbox{for $x\in L$,}
\end{eqnarray}
\begin{equation}\label{eq5.10}
(id\otimes \Delta)\Delta =(\Delta\otimes id)\Delta ,
\end{equation}
\begin{equation}\label{eq5.11}
\varepsilon (q)=1, \qquad  \varepsilon (x)=0 \qquad\mbox{for $x\in L$,}
\end{equation}
\begin{equation}\label{eq5.12}
(\varepsilon \otimes id)\Delta (a)=1\otimes \sigma (a)\qquad\mbox{for $a\in U^{6-th}(L)$,}
\end{equation}
and
\begin{equation}\label{eq5.13} 
(id\otimes \varepsilon)\Delta (a)=\sigma (a)\otimes 1\qquad\mbox{for $a\in U^{6-th}(L)$,}
\end{equation}
where $\sigma$ is the linear map defined in Proposition~\ref{pr5.2}.
\end{proposition}

\medskip
\noindent
{\bf Proof} Clearly, $\Big (U^{6-th}(L)\otimes U^{6-th}(L), q\otimes q\Big)$ is an invariant algebra induced by the idempotent $q\otimes q$. First, we prove that the linear map 
\begin{equation}\label{eq5.14} 
\Delta :L\to Lie\bigg(\Big (U^{6-th}(L)\otimes U^{6-th}(L), q\otimes q\Big), 
[\, ,\,]_{6, k}\bigg)
\end{equation}
defined by (\ref{eq5.9}) is a Lie algebra homomorphism. Let $x$ and $y$ be elements of $L$. Regarding $L$ as a subalgebra in 
$Lie\bigg(\Big (U^{6-th}(L)\otimes U^{6-th}(L), q\otimes q\Big), [\, ,\,]_{6, k}\bigg)$, we have
\begin{eqnarray}\label{eq5.16}
&&\Delta([x, y])\nonumber\\
&=&\Big([x, y]+kq[x, y]-[x, y]q\Big)\otimes 1
+1\otimes \Big([x, y]+kq[x, y]-[x, y]q\Big)+\nonumber\\
&&\qquad +(1-k)q[x, y]\otimes q+(1-k)q\otimes q[x, y]\nonumber\\
&=&(xy-yx-xyq+yxq+k^2qxy-k^2qyx)\otimes 1+\nonumber\\
&&\qquad +1\otimes (xy-yx-xyq+yxq+k^2qxy-k^2qyx)+\nonumber\\
&&\qquad +(k-k^2)(qxy-qyx)\otimes q+(k-k^2)q\otimes (qxy-qyx).
\end{eqnarray}

In the Lie algebra 
$Lie\bigg(\Big (U^{6-th}(L)\otimes U^{6-th}(L), q\otimes q\Big), [\, ,\,]_{6, k}\bigg)$, we have
\begin{eqnarray}\label{eq5.19}
&&[\Delta(x), \Delta(y)]_{6, k}=\Delta(x)\Delta(y)-\Delta(y)\Delta(x)+\nonumber\\
&&\qquad +(k-1)(q\otimes q)\Big(\Delta(x)\Delta(y)-\Delta(y)\Delta(x)\Big).
\end{eqnarray}

Using the following fact:
$$
(x+kqx-xq)(y+kqy-yq)=(xy+k^2qxy-xyq),
$$
we get
\begin{eqnarray}\label{eq5.20}
&&\Delta(x)\Delta(y)\nonumber\\
&=&(xy+k^2qxy-xyq)\otimes 1+\underbrace{(x+kqx-xq)\otimes (y+kqy-yq)}_{\star}
+\nonumber\\
&&+(1-k^2)qxy\otimes q+(1-k^2)q\otimes qxy+\nonumber\\
&&+\underbrace{(y+kqy-yq)\otimes (x+kqx-xq)}_{\star}+1\otimes (xy+k^2qxy-xyq)\nonumber\\
&&+\underbrace{(1-k^2)qy\otimes qx}_{\star}+\underbrace{(1-k^2)qx\otimes qy}_{\star}
\end{eqnarray}

Note that the sum of the terms with the label $\star$ in (\ref{eq5.20}) is symmetric in $x$ and $y$. Hence, the terms with the label $\star$ disappear in the difference 
$\Delta(x)\Delta(y)-\Delta(y)\Delta(x)$. Thus, (\ref{eq5.20}) implies that
\begin{eqnarray}\label{eq5.21}
&&\Delta(x)\Delta(y)-\Delta(y)\Delta(x)\nonumber\\
&=&(xy+k^2qxy-xyq-yx-k^2qyx+yxq)\otimes 1+\nonumber\\
&&+(1-k^2)(qxy-qyx)\otimes q+(1-k^2)q\otimes (qxy-qyx)+\nonumber\\
&&+1\otimes (xy+k^2qxy-xyq-yx-k^2qyx+yxq).
\end{eqnarray}

By (\ref{eq5.21}), we get
\begin{eqnarray}\label{eq5.22}
&&(q\otimes q)\Big(\Delta(x)\Delta(y)-\Delta(y)\Delta(x)\Big)\nonumber\\
&=&(qxy-qyx)\otimes q+q\otimes (qxy-qyx).
\end{eqnarray}

Using  (\ref{eq5.21}) and (\ref{eq5.22}), (\ref{eq5.19}) becomes
\begin{eqnarray}\label{eq5.23}
&&[\Delta(x), \Delta(y)]_{6, k}\nonumber\\
&=&(xy+k^2qxy-xyq-yx-k^2qyx+yxq)\otimes 1+\nonumber\\
&& +(k-k^2)(qxy-qyx)\otimes q+(k-k^2)q\otimes (qxy-qyx)+\nonumber\\
&& +1\otimes (xy+k^2qxy-xyq-yx-k^2qyx+yxq).
\end{eqnarray}

It follows from (\ref{eq5.16}) and (\ref{eq5.23}) that
\begin{equation}\label{eq5.24} 
\Delta ([x, y])=[\Delta(x), \Delta(y)]_{6, k} \quad\mbox{for $x, y\in L$,}
\end{equation}
which proves that the linear map $\Delta$ given by (\ref{eq5.9}) is a Lie algebra homomorphism. By the property of the enveloping$^{6-th}$ algebra $U^{6-th}(L)$, the linear map $\Delta$ given by (\ref{eq5.9}) can be extended to an invariant homomorphism
$$
\Delta : \Big(U^{6-th}(L), q\Big)\to \Big(U^{6-th}(L)\otimes U^{6-th}(L), q\otimes q\Big).
$$
Thus, $\Delta : U^{6-th}(L)\to U^{6-th}(L)\otimes U^{6-th}(L)$ is an algebra homomorphism such that both (\ref{eq5.8}) and (\ref{eq5.9}) hold.

\bigskip
Next, we prove (\ref{eq5.10}) which means that $\Delta$ is a comultiplication. Clearly, we have
\begin{equation}\label{eq5.25} 
(id\otimes \Delta)\Delta (1)=1\otimes 1\otimes 1=(\Delta\otimes id)\Delta (1)
\end{equation}
and
\begin{equation}\label{eq5.26} 
(id\otimes \Delta)\Delta (q)=q\otimes q\otimes q=(\Delta\otimes id)\Delta (q).
\end{equation}

For $x\in L$, we have
\begin{eqnarray}\label{eq5.27} 
&&(id\otimes \Delta)\Delta (x)\nonumber\\
&=&(x+kqx-xq)\otimes \Delta (1)+
1\otimes \Big(\Delta (x)+k\Delta (qx)-\Delta (xq)\Big)+\nonumber\\
&& +(1-k)qx\otimes \Delta (q)+(1-k)q\otimes \Delta (qx).
\end{eqnarray}

Since
\begin{eqnarray*}\label{eq5.28} 
&&\Delta (qx)=(q\otimes q)\Delta (x)=\Delta (x)(q\otimes q)\nonumber\\
&=&\Delta (xq)=qx\otimes q+q\otimes qx 
\end{eqnarray*}
and
\begin{eqnarray*}\label{eq5.29} 
&&\Delta (x)+k\Delta (qx)-\Delta (xq)=\Delta (x)+(k-1)\Delta (qx)\nonumber\\
&=&(x+kqx-xq)\otimes 1+1\otimes (x+kqx-xq),
\end{eqnarray*}
(\ref{eq5.27}) becomes
\begin{eqnarray}\label{eq5.30} 
&&(id\otimes \Delta)\Delta (x)\nonumber\\
&=&(x+kqx-xq)\otimes 1\otimes 1+1\otimes (x+kqx-xq)\otimes 1+\nonumber\\
&&+1\otimes 1\otimes (x+kqx-xq)+(1-k)qx\otimes q\otimes q+\nonumber\\
&&+(1-k)q\otimes qx\otimes q+(1-k)q\otimes q\otimes qx.
\end{eqnarray}

Similarly, we have
\begin{eqnarray}\label{eq5.31} 
&&(\Delta\otimes id)\Delta (x)\nonumber\\
&=&
\Big(\Delta (x)+k\Delta (qx)-\Delta (xq)\Big)\otimes 1+\Delta (1)\otimes (x+kqx-xq)+\nonumber\\
&&+(1-k)\Delta (qx)\otimes q+(1-k)\Delta (q)\otimes qx\nonumber\\
&=&(x+kqx-xq)\otimes 1\otimes 1+1\otimes (x+kqx-xq)\otimes 1+\nonumber\\
&&+1\otimes 1\otimes (x+kqx-xq)+(1-k)qx\otimes q\otimes q+\nonumber\\
&&+(1-k)q\otimes qx\otimes q+(1-k)q\otimes q\otimes qx.
\end{eqnarray}

By (\ref{eq5.30}) and (\ref{eq5.31}), we get
\begin{equation}\label{eq5.32} 
(id\otimes \Delta)\Delta (x)=(\Delta\otimes id)\Delta (x)\quad\mbox{for $x\in L$.}
\end{equation}

Since the algebra $U^{6-th}(L)$ is generated by the set $\{1, q\}\cup L$, (\ref{eq5.10}) holds by 
(\ref{eq5.25}), (\ref{eq5.26}) and (\ref{eq5.32}).

\bigskip
We now construct the algebra homomorphism $\varepsilon :  U^{6-th}(L)\to\mathbf{k}$. Note that the identity $1$ of the field $\mathbf{k}$ is an idempotent, and $(\mathbf{k}, 1)$ is an invariant algebra induced by the idempotent $1$. If 
$\varepsilon :  L\to Lie \Big((\mathbf{k}, 1), [ \, , \, ]_{6, k}\Big)$ is the linear map defined by
$$
\varepsilon (x):=0 \qquad\mbox{for $x\in L$,}
$$
then $\varepsilon$ clearly is a Lie algebra homomorphism. Hence, $\varepsilon$ can be extended to an invariant homomorphism $\varepsilon : \Big(U^{6-th}(L), q\Big)\to (\mathbf{k}, 1)$. Thus, (\ref{eq5.11}) holds.

\medskip
Finally, both (\ref{eq5.12}) and (\ref{eq5.13}) follow from the definitions of $\varepsilon$ and 
$\Delta$.

\hfill\raisebox{1mm}{\framebox[2mm]{}}

\medskip
If $V$ and $W$ are vector spaces over a field $\mathbf{k}$, then the {\bf twist map} 
$$\tau : V\otimes W\to W\otimes V$$ 
is defined by $\tau (v\otimes w):=w\otimes v$, where $v\in V$ and $w\in W$.

\medskip
First, we introduce  the concept of a bialgebra with $\sigma$-counit.  

\begin{definition}\label{def5.3} A vector space $H$ over a field $\mathbf{k}$ is called a 
{\bf bialgebra with $\sigma$-counit} if there exist five linear maps $m: H\otimes H\to H$, 
$u: \mathbf{k}\to H$, $\Delta : H\to H\otimes H$,
$\varepsilon : H\to \mathbf{k}$ and $\sigma : H\to H$  such that the following eight diagrams are commutative:
\begin{itemize} 
\item associativity
$$
\begin{array}{ccc}
H\otimes H\otimes H
&\stackrel{m\otimes id}{\longrightarrow}&H\otimes H\\
id\otimes m\Bigg\downarrow&&\Bigg\downarrow m\\
H\otimes H&\stackrel{m}{\longrightarrow}&H
\end{array}
$$
\item unit
$$
\begin{array}{ccccc}
\mathbf{k}\otimes H
&\stackrel{u\otimes id}{\longrightarrow}&H\otimes H&
\stackrel{id\otimes u}{\longleftarrow}&H\otimes \mathbf{k}\\
\Bigg\downarrow&&\Bigg\downarrow m&&\Bigg\downarrow\\
H&\stackrel{id}{=}&H&\stackrel{id}{=}&H
\end{array}
$$
\item coassociativity
$$
\begin{array}{ccc}
H&\stackrel{\Delta}{\longrightarrow}&H\otimes H\\
\Delta\Bigg\downarrow&&\Bigg\downarrow id\otimes\Delta\\
H\otimes H&\stackrel{\Delta\otimes id}{\longrightarrow}&H\otimes H\otimes H
\end{array}
$$
\item $\sigma$-counit
$$
\begin{array}{ccccc}
\mathbf{k}\otimes H&\longleftarrow&H&\longrightarrow&H\otimes \mathbf{k}\\
\varepsilon\otimes id\Bigg\uparrow&&\Bigg\uparrow\sigma&&\Bigg\uparrow id\otimes \varepsilon\\
H\otimes H&\stackrel{\Delta}{\longleftarrow}&H&\stackrel{\Delta}{\longrightarrow}
&H\otimes H
\end{array}
$$
\item $\Delta$ is an algebra homomorphism
$$
\begin{array}{ccc}
H\otimes H&\stackrel{m}{\longrightarrow}H\stackrel{\Delta}{\longrightarrow}&H\otimes H\\
\Delta\otimes \Delta\Bigg\downarrow&&\Bigg\uparrow m\otimes m\\
H\otimes H\otimes H\otimes H&\stackrel{id\otimes \tau\otimes id}{\longrightarrow}&H\otimes H\otimes H\otimes H
\end{array}
$$
$$
\begin{array}{ccc}
H&\stackrel{\Delta}{\longrightarrow}&H\otimes H\\
u\Bigg\uparrow&&\Bigg\uparrow u\otimes u\\
\mathbf{k}&\longrightarrow&\mathbf{k}\otimes\mathbf{k}
\end{array}
$$
\item $\varepsilon$ is an algebra homomorphism
$$
\begin{array}{cc}
\begin{array}{ccc}
H\otimes H&\stackrel{\varepsilon\otimes\varepsilon}{\longrightarrow}&\mathbf{k}\otimes\mathbf{k}\\
m\Bigg\downarrow&&\Bigg\downarrow \\
H&\stackrel{\varepsilon}{\longrightarrow}&\mathbf{k}
\end{array}\qquad\qquad&
\begin{array}{ccc}
H&=&H\\
u\Bigg\uparrow&&\Bigg\downarrow \varepsilon\\
\mathbf{k}&=&\mathbf{k}
\end{array}
\end{array}
$$
\end{itemize}
The linear map $\varepsilon$ is called the {\bf $\sigma$-counit}. A bialgebra H with 
$\sigma$-counit is also denoted by $H{m, u, \Delta \choose  \varepsilon , \sigma }$ if the five linear maps have to be indicated explicitly.
\end{definition}

Based on bialgebras with $\sigma$-counit, we now introduce Hopf-like algebras in the following 

\begin{definition}\label{def4.2} A  bialgebra $H=H{m, u, \Delta \choose  \varepsilon , \sigma }$  with $\sigma$-counit over a field $\mathbf{k}$ is called a 
{\bf Hopf-like algebra} if there exists a linear map $S : H\to H$ such that the following diagram is commutative:
\begin{equation}\label{eq588}
\begin{array}{ccccccc}
H&\stackrel{\Delta}{\longrightarrow}&H\otimes H&\stackrel{S\otimes id}{\longrightarrow}
&H\otimes H&\stackrel{m}{\longrightarrow}&H\\
\Big|\Big|&&&&&&\Big|\Big|\\
H&\stackrel{S}{\longrightarrow}&H&\stackrel{\varepsilon}{\longrightarrow}&\mathbf{k}
&\stackrel{u}{\longrightarrow}&H\\
\Big|\Big|&&&&&&\Big|\Big|\\
H&\stackrel{\Delta}{\longrightarrow}&H\otimes H&\stackrel{id\otimes S}{\longrightarrow}
&H\otimes H&\stackrel{m}{\longrightarrow}&H
\end{array}
\end{equation}
\end{definition}
The linear map $S$ is called the {\bf antipode -like} of the Hopf-like  algebra $H$.

\medskip
The Hopf algebra structure of the ordinary enveloping algebra for a Lie algebra is generalized in the following

\begin{proposition}\label{pr5.4} The enveloping$^{6-th}$ algebra $U^{6-th}(L)$ of a Lie algebra $L$ is a Hopf-like algebra, where the comultiplication $\Delta$, the $\sigma$-counit 
$\varepsilon$ and the antipode-like $S$ are given by Proposition~\ref{pr5.1} and Proposition~\ref{pr5.2}.
\end{proposition}

\medskip
\noindent
{\bf Proof} It follows from Proposition~\ref{pr5.1} and Proposition~\ref{pr5.2} that we need only to prove that the diagram in (\ref{eq588}) is commutative, i.e.
\begin{equation}\label{eq5.38}
m(S\otimes id)\Delta =u\varepsilon S=m(id\otimes S)\Delta .
\end{equation}

\medskip
Let {\bf 1} be the identity of $U^{6-th}(L)$. Clearly, we have
\begin{equation}\label{eq5.39}
m(S\otimes id)\Delta ({\bf 1})={\bf 1}=u\varepsilon S({\bf 1})
=m(id\otimes S)\Delta ({\bf 1}) 
\end{equation}
and
\begin{equation}\label{eq5.40}
m(S\otimes id)\Delta (q)=u\varepsilon S(q)=m(id\otimes S)\Delta (q). 
\end{equation}

\medskip
For $x\in L$, we have
\begin{eqnarray}\label{eq5.41}
&&m(S\otimes id)\Delta (x)\nonumber\\
&=&m(S\otimes id)\Big\{(x+kqx-xq)\otimes 1+1\otimes (x+kqx-xq)+\nonumber\\
&&\quad  +(1-k)qx\otimes q+(1-k)q\otimes qx\Big\}\nonumber\\
&=&S(x)+kS(x)S(q)-S(q)S(x)+(x+kqx-xq)+\nonumber\\
&&\quad +(1-k)S(x)\underbrace{S(q)q}_{0}+(1-k)\underbrace{S(q)q}_{0}x\nonumber\\
&=&\underbrace{-\frac1k x-kqx+\frac1k xq}
+k\left(-\frac1k x-kqx+\frac1k xq\right)(1-q)+\nonumber\\
&&\quad -\underbrace{(1-q)\left(-\frac1k x-kqx+\frac1k xq\right)}+x+kqx-xq\nonumber\\
&=&q\left(-\frac1k x-kqx+\frac1k xq\right)+(-x-k^2qx+xq)(1-q)+x+kqx-xq\nonumber\\
&=&-kqx-x(1-q)-k^2qx(1-q)+xq(1-q)+x+kqx-xq\nonumber\\
&=&0,
\end{eqnarray}
\begin{equation}\label{eq5.42}
u\varepsilon S(x)=u\varepsilon\left(-\frac1k x-kqx+\frac1k xq\right)=u(0)=0
\end{equation}
and
\begin{eqnarray}\label{eq5.43}
&&m(id\otimes S)\Delta (x)\nonumber\\
&=&m(id\otimes S)\Big\{(x+kqx-xq)\otimes 1+1\otimes (x+kqx-xq)+\nonumber\\
&&\quad  +(1-k)qx\otimes q+(1-k)q\otimes qx\Big\}\nonumber\\
&=&x+kqx-xq+S(x)+kS(x)S(q)-S(q)S(x)+\nonumber\\
&&\quad +(1-k)\underbrace{qxS(q)}_{0}+(1-k)\underbrace{qS(qx)}_{0}\nonumber\\
&=&0.
\end{eqnarray}

By (\ref{eq5.41}), (\ref{eq5.42}) and (\ref{eq5.43}), we get
\begin{equation}\label{eq5.44}
m(S\otimes id)\Delta (x)=u\varepsilon S(x)=m(id\otimes S)\Delta (x)\qquad\mbox{for $x\in L$.}
\end{equation}

\medskip
Since the algebra $U^{6-th}(L)$ is generated by the set $\{1, q\}\cup L$, (\ref{eq5.38}) holds by 
(\ref{eq5.39}), (\ref{eq5.40}) and (\ref{eq5.44}).

\hfill\raisebox{1mm}{\framebox[2mm]{}}

\end{document}